\documentclass{amsart}

\usepackage{amsfonts}
\usepackage{amsmath}
\usepackage{amsthm}
\usepackage{amssymb}
\usepackage[english]{babel}
\usepackage{graphics}
\usepackage{latexsym}
\usepackage{longtable} 
\usepackage{mathrsfs} 
\usepackage{graphicx}
\usepackage{wrapfig} %\usepackage{txfonts}
\usepackage[pdftex,colorlinks=true,linkcolor=blue,citecolor=blue]{hyperref}
\usepackage{enumitem}
\usepackage{esint}

\newtheorem{theorem}{Theorem}
\newtheorem{corollary}[theorem]{Corollary}
\newtheorem{lemma}[theorem]{Lemma}
\newtheorem{proposition}[theorem]{Proposition}

\newtheorem{claim}[theorem]{Claim}
\newtheorem{example}[theorem]{Example}

\theoremstyle{definition}
\newtheorem{definition}[theorem]{Definition}
\newtheorem{remark}[theorem]{Remark}

\newcommand{\mL}{\mathcal{L}}
\newcommand{\mH}{\mathcal{H}}
\newcommand{\mF}{\mathcal{F}}

\newcommand{\mW}{\mathcal{W}}
\newcommand{\mX}{\mathfrak{X}}

\newcommand{\mM}{\mathcal{M}}
\newcommand{\mO}{\mathcal{O}}

\newcommand{\E}{\mathrm{E}}
\newcommand{\D}{\mathrm{D}}
\renewcommand{\H}{\mathrm{H}}

\newcommand{\A}{\mathrm{A}}
\newcommand{\B}{\mathrm{B}}
\newcommand{\K}{\mathrm{K}}
\renewcommand{\L}{\mathrm{L}}
\renewcommand{\H}{\mathrm{H}}
\newcommand{\G}{\mathrm{G}}

\newcommand{\R}{\mathbb{R}}
\newcommand{\C}{\mathbb{C}}
\newcommand{\N}{\mathbb{N}}

\newcommand{\mB}{\mathbb{B}}

\renewcommand{\L}{\mathrm{L}}

\newcommand{\noi}{\noindent}
\newcommand{\ms}{\medskip}

\newcommand{\al}{\alpha}
\newcommand{\be}{\beta}
\newcommand{\ga}{\gamma}
\newcommand{\Ga}{\Gamma}
\newcommand{\de}{\delta}

\newcommand{\e}{\varepsilon}
\newcommand{\si}{\sigma}

\newcommand{\la}{\lambda}

\newcommand{\ka}{\kappa}
\newcommand{\Om}{\Omega}
\newcommand{\om}{\omega}

\newcommand{\av}{-\hspace{-10.5pt}\displaystyle\int}

\newcommand{\weak }{\, -\!\!\!\!\!-\!\!\!\!\rightharpoonup}
\newcommand{\weakstar }{ \overset{\, *_{\phantom{|}}}{{\smash{\, -\!\!\!\!-\!\!\!\!\rightharpoonup}}\, } }

\newcommand{\larrow}{\longrightarrow}
\newcommand{\ot}{\otimes}

\newcommand{\LL}{\text{\LARGE$\llcorner$}}
\newcommand{\p}{\partial}
\newcommand{\sub}{\subseteq}
\newcommand{\set}{\setminus}
\newcommand{\by}{\times}

\renewcommand{\div}{\mathrm{div}}

\newcommand{\bt}{\begin{theorem}}\newcommand{\et}{\end{theorem}}
\newcommand{\bd}{\begin{definition}}\newcommand{\ed}{\end{definition}}
\newcommand{\bl}{\begin{lemma}}\newcommand{\el}{\end{lemma}}
\newcommand{\beq}{\begin{equation}}\newcommand{\eeq}{\end{equation}}
\newcommand{\bc}{\begin{claim}}\newcommand{\ec}{\end{claim}}
\newcommand{\bex}{\begin{example}}\newcommand{\eex}{\end{example}}
\newcommand{\bcor}{\begin{corollary}}\newcommand{\ecor}{\end{corollary}}
\newcommand{\bp}{\begin{proof}}\newcommand{\ep}{\end{proof}}

\numberwithin{equation}{section}
\begin{document}

\title[An $L^\infty$ approach to inverse optical tomography]{Inverse optical tomography through PDE constrained optimisation in $L^\infty$}
 
\author{Nikos Katzourakis}

\address{Department of Mathematics and Statistics, University of Reading, Whiteknights, PO Box 220, Reading RG6 6AX, United Kingdom}

\email{n.katzourakis@reading.ac.uk}

  \thanks{\!\!\!\!\!\!\!\!\texttt{The author has been partially financially supported by the EPSRC grant EP/N017412/1}}
  
%\subjclass[2010]{35D99, 35D40, 35J47, 35J47, 35J92, 35J70, 35J99}

\date{}


\begin{abstract} Fluorescent Optical Tomography (FOT) is a new bio-medical imaging method with wider industrial applications. It is currently intensely researched since it is very precise and with no side effects for humans, as it uses non-ionising red and infrared light. Mathematically, FOT can be modelled as an inverse parameter identification problem, associated with a coupled elliptic system with Robin boundary conditions. Herein we utilise novel methods of Calculus of Variations in $L^\infty$ to lay the mathematical foundations of FOT which we pose as a PDE-constrained minimisation problem in $L^p$ and $L^\infty$.

\end{abstract}

\maketitle

\tableofcontents

\section{Introduction}   \label{Section1}

Fluorescent Optical Tomography (FOT) is a relatively new and still evolving bio-medical imaging method, which also has wider industrial applications. FOT is currently being very intensely studied, as it presents some advantages over more classical imaging methods which use X-rays, Gamma-rays, electromagnetic radiation and ultrasounds.

The goal of FOT is to reconstruct interior optical properties of an object (e.g.\ living tissue) by using light in the red (visible) and infrared (invisible) range. The principal current use of FOT is in medical applications (e.g.\ cancer tumours diagnosis and general prevention of various diseases), as well as in industrial applications (e.g.\ detecting structural flaws in superconductors), see e.g.\ \cite{ARH, AMM, Ar, BJ, BJ2, FP, FES, GHA, GEZ, JBS, ZRWN, ZG}.

FOT improves on quite a few shortcomings of other popular imaging methods. The vast majority of currently available bio-imaging techniques image merely tissue structure variations created by tumours. However, since some features imaged are not specific to the presence of actual tumour cells, the unavoidable imaging of secondary effects might lead to false diagnoses. Additionally, imaging methods using X-rays and Gamma-rays actually use ionising radiation, which is harmful for humans and animals as it is potentially cancer-inducing itself. On the other hand, FOT is an imaging method which does not use harmful radiation and can be made specific to the presence of designated cell types. Therefore, FOT is more precise and with no side effects for humans.

Technically, the aim of FOT is to reconstruct the fluorophore distribution in a solid body from measurements of light intensity through detectors placed on the boundary. The highly diffusive nature of light propagation implies that in fact FOT forms a highly nonlinear and severely ill-posed inverse problem, hence mathematically it is a very challenging problem. FOT can be modelled by a coupled system of PDEs (partial differential equations) with $\C$-valued solutions and coefficients. The goal is to reconstruct a space-varying parameter in the system of PDEs in the interior of a body (e.g.\ living tissue).

Mathematically, FOT can be modelled as follows. Let $\Om \sub \R^n$ be an open bounded set with $C^1$ boundary $\p\Om$ and $n\geq 3$. In medical applications $n=3$, but from the mathematical viewpoint we may include greater dimensions without any ramifications. 
 A fluorescent dye is injected into the body $\Om$ and in order to determine the dye concentration $\xi=\xi(x)$, the body is illuminated by a red light source $s=s(x)$ placed on the boundary $\p\Om$. The wavelength of the light is adjusted to the excitation wavelength of the dye, in order to force it to fluoresce. The light diffuses inside the body, and wherever dye is present, fluorescent light in the infrared range is emitted that can then be detected again at the body surface using a camera and appropriate infrared filters. The goal is then to reconstruct the distribution $\xi=\xi(x)$ of the dye, from these obtained surface images.

Specifically, for time-periodic light sources modulated at a specific frequency, the following system of PDEs describes at any $x\in \Om$ the $\C$-valued photon fluences $u=u(x)$ at the excitation wavelength and $v=v(x)$ at the fluorescent wavelength:
\beq \label{1.1}
\left\{\ \ 
\begin{array}{ll}
  -\div(\A_\xi \D u)\,+\, k_\xi u\, =\, S, & \ \ \text{ in }\Om, \ms
 \\
  -\div(\A_\xi \D v)\,+\, k_\xi v\, =\, \xi h u, & \ \ \text{ in }\Om,  \ms
\\
 \ \  (\A_\xi \D u) \cdot \mathrm n\,+\, \ga u\, =\, s, & \ \ \text{ on }\p\Om,  \ms
\\
 \ \ (\A_\xi \D v) \cdot \mathrm n\,+\, \ga v\, =\, 0, & \ \ \text{ on }\p\Om.
\end{array}
\right.
\eeq
Here $u,v : \Om \larrow \C$ are the solutions,  {$\D u, \D v$ are the gradients of the solutions (namely $\D = (\D_1,...,\D_n)$ and $\D_k \equiv \p/\p x_k$)} and $\mathrm n : \p\Om \larrow \R^n$ is the outer unit normal vector field on the boundary, whilst the $\xi$-dependent coefficients $\A_\xi,k_\xi$ and the coefficients $h,s,S,\xi,\ga$ satisfy $\ga>0$ and 
\beq \label{1.2}
\begin{split}
\A_\xi \ & : \ \Om \larrow \R^{n\by n}_+,
\\
k_\xi, h,S\ & : \ \Om \larrow \C,
\\
s\ & : \ \p\Om \larrow \C,
\\
\xi\ & : \ \Om \larrow [0,\infty) 
\end{split}
\eeq
with the real part of $k_\xi$ being positive. We note that our general PDE notation will be either self-explanatory, or otherwise standard, as e.g.\ in the textbooks \cite{E,KV}. For instance, $\R^{n\by n}_+$ symbolises the cone of real non-negative $n\by n$ matrices. In imaging applications where $n=3$, the coefficients above take the following form:
\beq
\label{1.1A}
\left\{
\ \ \ 
\begin{split}
 \A_\xi (x) \, &=\,  \Big(3\big(\mu_{ai}(x) \,+\, \mu'_{s}(x)  \,+\, \xi(x) \big)\Big)^{\!-1}\, \mathrm I_n,
\\
k_\xi(x)\,&=\, \mu_{ai}(x)\,+\, \xi(x)\,+\, {i\om}\mathrm c^{-1},
\\
h(x)\,&=\, {\phi}\big(1-i\om \tau(x)\big)^{-1},
\end{split}
\right.
\eeq
where $\mathrm I_n$ is the identity matrix in $\R^n$, the diffusion coefficient $\A_\xi$ describes the diffusion of photons, $\mu_{ai}$ is the absorption coefficient due to the endogenous chromophores, $\xi$ is the absorption coefficient due to the exogenous fluorophore, $\mu'_s$ is the reduced scattering coefficient, $\phi$ is the quantum efficiency of the fluorophore, $\tau$ is the fluorophore lifetime and $\om$ is the modulated light frequency and $\mathrm c$ the speed of light. Finally, $S,s$ are the light sources. In applications, some authors model the problem with either boundary sources or interior light sources (see e.g.\ \cite{FES} versus \cite{BJ, BJ2}). Mathematically, we may include both types of sources without difficulty.

To the best of our knowledge, although the FOT problem has been extensively studied computationally and numerically, it has not been been considered from the purely analytical viewpoint. In this paper we utilise novel methods of Calculus of Variations in $L^\infty$ in order to lay the rigorous mathematical foundations of the FOT problem. Motivated by developments underpinning the papers \cite{K, KM, KP1, KPa}, we pose FOT as a minimisation problem in $L^\infty$ with PDE constraints as well as unilateral constraints, studying the direct as well as the inverse FOT problem, both in $L^p$ for finite $p$ and in $L^\infty$. Further, we derive the relevant variational inequalities in $L^p$ for finite $p$ and in $L^\infty$ that the constrained minimisers satisfy, which involve (generalised) Lagrange multipliers. Additionally, we prove convergence of the corresponding $L^p$ structures to the limiting $L^\infty$ structures as $p\to\infty$, in a certain fashion that will become clear later. 

Calculus of Variations in $L^\infty$ is a modern area initiated by Aronsson in the 1960s (see \cite{A1}-\cite{A4}) who was the first to consider variational problems of functionals which are defined as a supremum. For a general pedagogical introduction we refer e.g.\ to \cite{C,K4}. Except for the endogenous mathematical interest, $L^\infty$ cost (or error) functionals are important for optimisation and control applications because by minimising the supremum of the cost rather than its average (as e.g.\ in the case of least square $L^2$ costs), we obtain better results as we achieve uniform smallness of the error and small area spike deviations are a priori excluded. Interesting theory and applications of $L^\infty$ variational problems can be found e.g.\ in \cite{AB, BBJ, BaJe, BJW1, BN, BP, CDP, GNP, MWZ, PP, P, RZ}.

For our purposes in this paper, it will be more convenient to rewrite \eqref{1.1} in vectorial rather than complex form because the real/vectorial form. Hence, for any complex valued function $f=f_R  + if_I: \Om \larrow \C$, we identify $f$ with $(f_R,f_I)^\top : \Om \larrow\R^2$ and we will consider the following generalisation of \eqref{1.1}-\eqref{1.1A}:
\beq \label{1.3}
\left\{\ \ 
\begin{array}{lll}
(a)\ \ & \ \  -\div(\A_\xi \overset{_\bullet}{\phantom{.}} \D u)\,+\, \K_\xi u\, =\, S, & \ \ \text{ in }\Om, \ms
\\
(b)&\ \     -\div(\A_\xi {\hspace{1pt}^{_\bullet}} \D v)\,+\, \K_\xi v\, =\, \xi \H u, & \ \ \text{ in }\Om,  \ms
\\
(c)&  \ \ \ \ (\A_\xi \hspace{1pt} \overset{_\bullet}{\phantom{.}} \D u)\hspace{1pt} \mathrm n\,+\, \ga u\, =\, s, & \ \ \text{ on }\p\Om,  \ms
\\
(d)& \ \ \ \  (\A_\xi \hspace{1pt} \overset{_\bullet}{\phantom{.}} \D v)\hspace{1pt} \mathrm n\,+\, \ga v\, =\, 0, & \ \ \text{ on }\p\Om.
\end{array}
\right.
\eeq
 {Note that although the complex form is seemingly more neat than the tangled vectorial form, it is nonetheless more appropriate in order to apply standard PDE methods and results, a priori estimates, maximum principle results, lower semicontinuity arguments, etc, which are all stated for real-valued or vector-valued solutions and PDE coefficients.} In \eqref{1.3} we have
\beq  \label{1.4}
\smash{\A_\xi \,:=\, \A + r(\cdot,\xi)\mathrm I_n, \ \ \ r(x,t)\,:=\, \frac{\la}{\kappa(x)+t},\ \ \ \A_\xi \overset{_\bullet}{\phantom{.}} \D u \,:=  \left[
\begin{array}{cc}
\D u_R \A_\xi
\\
\D u_I \A_\xi
\end{array}
\right],}
\eeq
with
\beq  \label{1.5}
\K_\xi \,:=\, \K \,+\, \xi \mathrm I_2, \ \ \ \K \, :=\,
\left[
\begin{array}{cc}
k_R  & -k_I
\\
k_I & k_R 
\end{array}
\right]
, \ \ \
\H\,:=\,
\left[
\begin{array}{cc}
h_R & -h_I
\\
h_I & h_R
\end{array}
\right],
\eeq
and 
\beq \label{1.6}
\left\{
\begin{split}
u,\, v,\, S\ & : \ \Om\larrow \R^2, \ \ \ \ \ \D u,\, \D v\ : \ \Om \larrow \R^{2\by n},
\\
\K,\,\H\ & : \ \Om \larrow \R^{2\by 2}, \ \ \ \ \ \ \ \ \ \, \A\ :\ \Om \larrow \R^{n\by n}_+,
\\
s\ & : \ \p \Om\larrow \R^2, \ \ \ \ \ \ \ \ \ \ \ \  \xi\  : \ \Om\larrow [0,\infty),
\\
& \ga,\la>0, \hspace{63pt} \ka \ :\ \Om \larrow (0,\infty).
\end{split}
\right.
\eeq
 {We note that when dealing with derivatives of vector-valued functions, the conventions about column/row vectors are important. Therefore, we interpret $u,v$ as column vectors and for each of the real scalar components $u_R,u_I,v_R,v_I$, we consider their gradients as row vectors. Similarly, for any matrix field $A = \sum_{k,l} A_{kl} e_k \ot e_l$ (where $\{e_1,...,e_n\}$ is the standard Euclidean basis of column vectors and $a \ot b := a b^\top$), it is customary to define
\[
  \div (A)\, :=\, \sum_{k,l} (\D_l A_{kl}) e_k,
\]
namely the divergence of a matrix is the column vector of divergences of its rows. Therefore, for $u$ we have
\[
\div \big( \A_\xi \overset{_\bullet}{\phantom{.}}\D u \big) \, =  \left[
\begin{array}{cc}
\div( \D u_R \A_\xi)
\\
\div( \D u_I \A_\xi)
\end{array}
\right] 
\, =\, 
\left[
\begin{array}{cc}
\sum_{l}\D_l \Big( \sum_k (\D_k u_R) (\A_\xi)_{kl} \Big)
\\
\sum_{l} \D_l \Big( \sum_k (\D_k u_I) (\A_\xi)_{kl} \Big)
\end{array}
\right] ,
\]
and similarly for $v$.}

The table of contents gives an idea of the organisation of the material in this paper. After this Introduction, in Section \ref{Section2} we delve into the study of the building stones of the FOT problem, namely of linear elliptic divergence PDE systems with Robin boundary conditions. To this end, we establish well-posedness of the direct problem for these in $L^p$ for all $p\in[2,\infty)$ (Theorems \ref{th1}-\ref{th2}). In Section \ref{Section3} we establish the well-posedness of the direct FOT problem in the appropriate $L^p$  spaces for all $p\in[2,\infty)$ (Theorem \ref{th3}). In Section \ref{Section4} we begin the study of the inverse FOT problem as a constrained minimisation problem with PDE and unilateral constraints, proving the existence of minimisers in $L^p$ for all $p\in[2,\infty]$ and the convergence of the $L^p$ minimisers to the $L^\infty$ minimiser as $p\to \infty$ (Theorem \ref{th4}).  In Section \ref{Section5} we establish the existence of generalised Lagrange multipliers to the $L^p$ constrained minimisation problem and the relevant variational inequalities, by invoking the infinite-dimensional counterpart of the Kuhn-Tucker theory (Theorem \ref{th7}). Finally, in Section \ref{Section6} we establish the corresponding results for the extreme case of variational inequalities in $L^\infty$ (Theorem \ref{th14}). Let us conclude with the note that in the most recent paper \cite{K2} an alternative approach has been followed to optical tomography which is mathematically distinct but of less relevance to applications, since the diffusion coefficients are assumed not to depend on the parameter $\xi$.

\ms

\section{Linear elliptic systems with Robin boundary conditions}  \label{Section2}

In this section we begin with an auxiliary result of independent interest, namely the well-posedness of general linear divergence systems with Robin boundary conditions. Below we start with the case of the $L^2$ theory, which effectively is an application of the Lax Milgram theorem (see e.g.\ \cite{E}).

\begin{theorem}[Well-posedness in $W^{1,2}$] \label{th1} Let $\Om \Subset \R^n$ be a domain with $C^1$ boundary and let also $\mathrm n : \p\Om \larrow \R^n$ be the outer unit normal vector field. Consider the next boundary value problem with Robin boundary conditions
\beq \label{2.1}
\left\{\ \ 
\begin{array}{ll}
  -\div(\B \hspace{1pt} \overset{_\bullet}{\phantom{.}} \D u)\,+\, \L u\, =\, f-\div F, & \ \ \text{ in }\Om, \ms
\\
\ \ (\B \hspace{1pt} \overset{_\bullet}{\phantom{.}} \D u - F)\hspace{1pt} \mathrm n\,+\, \ga u\, =\, g, & \ \ \text{ on }\p\Om,
\end{array}
\right.
\eeq
where the coefficients satisfy $\B : \Om \larrow \R^{n\by n}_+$, $\L : \Om \larrow \R^{2\by 2}$, $f :\Om \larrow \R^2$, $F : \Om \larrow \R^{2\by n}$, $g : \p \Om \larrow \R^2$ and $\ga>0$. We suppose that there exists $\be_0>0$ such that
\beq  \label{2.2}
\left\{
\ \ 
\begin{split}
& \B \in L^\infty(\Om;\R^{n\by n}_+),\ \ \ \B  \overset{_\bullet}{\phantom{.}} \D u \,:=  \left[
\begin{array}{cc}
\D u_R \B
\\
\D u_I \B
\end{array}
\right],\ \ \si(B) \sub \Big[\be_0,\frac{1}{\be_0}\Big],
\\ 
&\L \in L^\infty(\Om;\R^{2\by 2}),\ \ \ \ \L\,:=\,
\left[
\begin{array}{cc}
l_R & -l_I
\\
l_I & l_R
\end{array}
\right]
,\ \ \ l_R\geq \be_0
\end{split}
\right.
\eeq
and also that
\beq \label{2.3}
f \in L^2(\Om;\R^2),\ \ \ F \in L^2(\Om;\R^{2\by n}),\ \ \ g \in L^2(\p\Om;\R^2).
\eeq
Then, \eqref{2.1} has a unique weak solution in $W^{1,2}(\Om;\R^2)$ satisfying
\beq \label{2.4}
\left\{\ \ 
\begin{split}
\int_\Om\Big[ \B :(\D u^\top \D \psi) \, +\,  (\L u)\cdot \psi \Big] \, \mathrm{d}\mL^n \,+\, \int_{\p\Om}\big[ \ga u\cdot \psi\big] \, \mathrm{d}\mH^{n-1}
\\
= \int_\Om\Big[ f\cdot \psi \, +\,  F:\D \psi\Big] \, \mathrm{d}\mL^n\, +\, \int_{\p\Om}\big[ g \cdot \psi\big] \, \mathrm{d}\mH^{n-1},
\end{split}
\right.
\eeq   
for all $\psi \in W^{1,2}(\Om;\R^2)$. In addition, there exists $C>0$ depending only on the coefficients  and the domain such that
\beq 
\label{2.5}
\| u \|_{W^{1,2}(\Om)} \, \leq\, C\Big(\| f \|_{L^2(\Om)} +\, \| F \|_{L^2(\Om)} +\,\| g \|_{L^{2}(\p\Om)}\Big).
\eeq
\end{theorem}

In \eqref{2.4}, the notation ``$:$" symbolises the Euclidean (Frobenius) inner product in the matrix space $\R^{n\by n}$ and  ``$\cdot$"  the Euclidean inner product in $\R^{2}$.

\bp As we have already mentioned, the aim is to apply of the Lax Milgram theorem. (Note that the matrix $\L$ involved in the zeroth order term is not symmetric, therefore this is not a direct consequence of the Riesz representation theorem.) To this end, we define the bilinear form 
\[
\mB \ :\ \ W^{1,2}(\Om;\R^2) \by W^{1,2}(\Om;\R^2) \larrow \R
\]
by setting
\[
\mB[u,\psi] \,:=\, \int_\Om\Big[ \B :(\D u^\top \D \psi) \, +\,  (\L u)\cdot \psi \Big] \, \mathrm{d}\mL^n \,+\, \int_{\p\Om}\big[ \ga u\cdot \psi\big] \, \mathrm{d}\mH^{n-1}.
\]
Since $\B,\L$ are $L^\infty$, we immediately have by H\"older inequality that
\[
\big| \mB[u,\psi]  \big|\,\leq\, C \| u \|_{W^{1,2}(\Om)} \| \psi \|_{W^{1,2}(\Om)} 
\]
for some $C>0$ and all $u,\psi \in W^{1,2}(\Om;\R^2)$. Further, since
\[
\begin{split}
(\L u)\cdot u\, &=\, [u_R,\, u_I] \left[
\begin{array}{cc}
l_R & -l_I
\\
l_I & l_R
\end{array}
\right] \left[
\begin{array}{c}
u_R
\\
u_I
\end{array}
\right]
\\
& = \, l_R(u_R)^2 +\, l_R(u_I)^2
\\
&=\, l_R |u|^2
\\
&\geq \be_0 |u|^2,
\end{split}
\]
we have that 
\[
\begin{split}
\mB[u,u] \,& \geq\, \be_0\Big(\| \D u \|^2_{L^2(\Om)} +\, \| u \|^2_{L^2(\Om)}\Big) +\,\ga \| u \|^2_{L^{2}(\p\Om)} 
\\
& \geq\, \be_0\| u \|^2_{W^{1,2}(\Om)},
\end{split}
\]
for any $u \in W^{1,2}(\Om;\R^2)$. Hence, the bilinear form $\mB$ is bi-continuous and coercive, therefore the conditions of the Lax-Milgram theorem are satisfied. Thus, for any functional $\Phi \in (W^{1,2}(\Om;\R^2))^*$, there exists a unique $u \in W^{1,2}(\Om;\R^2)$ such that
\[
\ \ \ \mB[u,\psi] = \langle \Phi, \psi \rangle, \ \ \ \forall\ \psi \in W^{1,2}(\Om;\R^2).
\]
To conclude, it suffices to show that
\[
\langle \Phi, \psi \rangle\, := \int_\Om\Big[ f\cdot \psi \, +\,  F:\D \psi\Big] \, \mathrm{d}\mL^n\, +\, \int_{\p\Om}\big[ g \cdot \psi\big] \, \mathrm{d}\mH^{n-1}
\]
indeed defines an element of the dual space $(W^{1,2}(\Om;\R^2))^*$, and also establish the a priori estimate. Indeed, by the trace theorem in $W^{1,2}(\Om;\R^2)$, there is a $C>0$ depending on $\Om$ which allows to estimate
\[
\begin{split}
\big| \langle \Phi, \psi \rangle \big| \, & \leq\, \Big(\| f \|_{L^2(\Om)} +\, \| F \|_{L^2(\Om)} \Big) \| \psi \|_{W^{1,2}(\Om)} +\, \| g \|_{L^2(\p\Om)} \| \psi \|_{L^2(\p\Om)} 
\\
& \leq\, C\Big(\| f \|_{L^2(\Om)} +\, \| F \|_{L^2(\Om)} +\, \| g \|_{L^2(\p\Om)}\Big)  \| \psi \|_{W^{1,2}(\Om)},
\end{split}
\]
which shows that $\Phi$ is indeed a bounded linear functional. The choice of $\psi:=u$ together with Young inequality with $\e>0$, gives
\[
\begin{split}
\big| \langle \Phi, u \rangle \big| \, & \leq\, \frac{C^2}{4\e}\Big(\| f \|_{L^2(\Om)} +\, \| F \|_{L^2(\Om)} +\, \| g \|_{L^2(\p\Om)}\Big)^2 +\, \e \| u \|^2_{W^{1,2}(\Om)}.
\end{split}
\]
When combining the above estimate with the lower bound for $\mB[u,u]$, we conclude by choosing $\e<\be_0/2$.
\ep

Note that in the above proof we have employed the common practice of denoting by $C$ a generic constant whose value might change from step to step in an estimate. This practice will be utilised in the sequel freely. Now we show that the obtained unique weak solution to \eqref{2.1} is in fact more regular if the coefficients permit it.

\begin{theorem}[Well-posedness in $W^{1,p}$] \label{th2} In the setting of Theorem \ref{th1}, consider again the boundary value problem \eqref{2.1} with Robin boundary conditions. In addition to the assumptions in Theorem \ref{th1}, we also suppose that
\[
\B \in VMO(\R^n;\R^{n\by n}_+),\ \ f \in L^{\frac{np}{n+p}}(\Om;\R^2),\ \ F \in L^p(\Om;\R^{2\by n}),\ \ g \in L^p(\p\Om;\R^2),
\]
for some 
\[
p \, >\, \frac{2n}{n-2}.
\]
Then, the unique weak solution of \eqref{2.1} lies in the space $W^{1,p}(\Om;\R^2)$. In addition, there exists $C>0$ depending only on the coefficients, the domain and $p$ such that
\beq \label{2.7}
\| u \|_{W^{1,p}(\Om)} \, \leq\, C\Big(\| f \|_{L^{\frac{np}{n+p}}(\Om)} +\, \| F \|_{L^p(\Om)} +\,\| g \|_{L^p(\p\Om)}\Big).
\eeq
\end{theorem}

 {Note that the space $VMO(\R^n)$ is the standard space of (scalar) functions of vanishing mean oscillation, defined as the closure of the space of continuous functions $C^0_0(\R^n)$ vanishing at infinity in the space $BMO(\R^n)$ of bounded mean oscillation (see e.g. \cite{E}) and $VMO(\R^n;\R^{n\by n}_+)$ is the convex set of those positive matrix-valued maps with VMO components. The stronger assumption of the leading matrix coefficient being VMO instead of merely $L^\infty$ is required in order to obtain the higher integrability result stated above, which in turn relies on an estimate which is known to be true only under this stronger assumption.}

\bp The key ingredient is to apply the well-know estimate for the Robin boundary value problem for linear divergence elliptic equations which has the exact same form as \eqref{2.7}, but applies to the scalar version of the problem \eqref{2.1} with $\L\equiv0$, see \cite{ACGG, D, DK, Ge, Gh, KLS, N2}. Hence, we need some arguments to show that it is still true in the more general case of \eqref{2.1}. To this end, we rewrite \eqref{2.1} component-wise as
\[
\ \ \left\{ \ \ 
\begin{array}{ll}
  -\div(\D u_R \B) \, =\, \Big\{ f_R -\big(l_R u_R - l_Iu_I \big)\Big\}-\div F_R, & \ \ \text{ in }\Om, \ms
\\
\ \ (\D u_R \B - F_R)\cdot \mathrm n\,+\, \ga u_R\, =\, g_R, & \ \ \text{ on }\p\Om,
\end{array}
\right.
\]
\[
\left\{ \ \ 
\begin{array}{ll}
-\div(\D u_I \B) \, =\, \Big\{ f_I -\big(l_Ru_I + l_Iu_R \big)\Big\}-\div F_I, & \ \ \text{ in }\Om, \ms
\\
\ \ (\D u_I \B - F_I) \cdot \mathrm n\,+\, \ga u_I\, =\, g_I, & \ \ \text{ on }\p\Om.
\end{array}
\right.
\]
By applying the scalar estimate to the each of the boundary value problems separately, we have
\beq
\label{2.8}
\begin{split}
\| u_R \|_{W^{1,p}(\Om)} \, \leq\, & C\Big(\| f_R \|_{L^{\frac{np}{n+p}}(\Om)} +\, \| F_R \|_{L^p(\Om)} 
\\
&\ \ \ +\,\| g_R \|_{L^p(\p\Om)}  +\, \|\L \|_{L^\infty(\Om)} \| u \|_{L^{\frac{np}{n+p}}(\Om)}\Big)
\end{split}
\eeq
and also
\beq
\label{2.9}
\begin{split}
\| u_I \|_{W^{1,p}(\Om)} \,  \leq\, & C\Big(\| f_I \|_{L^{\frac{np}{n+p}}(\Om)} +\, \| F_I \|_{L^p(\Om)} 
\\
&\ \ \ +\,\| g_I \|_{L^p(\p\Om)}  +\, \|\L \|_{L^\infty(\Om)} \| u \|_{L^{\frac{np}{n+p}}(\Om)}\Big).
\end{split}
\eeq
Note now that since by assumption $p >\frac{2n}{n-2}$, we have
\[
2\, <\, \frac{np}{n+p}\,<\, p.
\]
Hence, by the $L^p$ interpolation inequalities, we can estimate
\[
\  \ \ \| u \|_{L^{\frac{np}{n+p}}(\Om)} \leq\,  \| u \|^\la_{L^{2}(\Om)} \,  \| u \|^{1-\la}_{L^{p}(\Om)}, \ \ \ \text{ where }\ \la = \frac{2p}{n(p-2)}.
\]
By the Young inequality (for $a,b\geq 0$, $\e>0$, $r>1$ and $r/(r-1)=r'$)
\beq
\label{2.10}
ab\, \leq \, \left\{\frac{r-1}{r}(\e r)^{\frac{1}{1-r}}\right\}b^{\frac{r}{r-1}}  \, +\, \, \e a^r,
\eeq
for the choice $r:=(1-\la)^{-1}$, we have
\[
r=\frac{n(p-2)}{p(n-2)-2n}\ ,\ \ \ \frac{r}{r-1}=\frac{n(p-2)}{2p} \ , \ \ \ 1-\la= \frac{n(p-2)}{p(n-2)-2n},
\]
and hence we can further estimate 
\begin{align}
 \| u \|_{L^{\frac{np}{n+p}}(\Om)} \, & \leq\,  \| u \|^\la_{L^{2}(\Om)} \,  \| u \|^{1-\la}_{L^{p}(\Om)}
 \nonumber\\
& =\,  \left(\| u \|_{L^{2}(\Om)}\right)^{\frac{2p}{n(p-2)}}   \left(\| u \|_{L^{p}(\Om)}\right)^{\frac{p(n-2)-2n}{n(p-2)}} 
\nonumber\\
& \leq\, \left\{\frac{r-1}{r}(\e r)^{\frac{1}{1-r}}\right\}  \! \Big( \!\left(\| u \|_{L^{2}(\Om)}\right)^{\frac{2p}{n(p-2)}}\! \Big)^{\!\frac{r}{r-1}} + \Big(\! \left(\| u \|_{L^{p}(\Om)}\right)^{\frac{p(n-2)-2n}{n(p-2)}} \! \Big)^{\!r}
\label{2.11}
 \\
 & =\,  \left\{\frac{2p}{n(p-2)}\left( \frac{\e n (p-2)}{p(n-2)-2n} \right)^{ \!\! -\frac{p(n-2)-2n}{2p} }\right\} \| u \|_{L^{2}(\Om)} \,+\, \e \| u \|_{L^{p}(\Om)} 
\nonumber \\
  & =: \,  C(\e,p,n) \| u \|_{L^{2}(\Om)} +\, \e \| u \|_{L^{p}(\Om)}. \phantom{\Big|} \nonumber
  \end{align}
By \eqref{2.8}, \eqref{2.9} and \eqref{2.11}, by choosing $\e>0$ small enough, we infer that 
\[
\| u \|_{W^{1,p}(\Om)} \, \leq\, C\Big(\| f \|_{L^{\frac{np}{n+p}}(\Om)} + \| F \|_{L^p(\Om)} +\,\| g \|_{L^p(\p\Om)} +\,\| u \|_{L^2(\Om)}\Big).
\]
The desired estimate \eqref{2.7} ensues by combining the above estimate with our earlier $W^{1,2}$ estimate \eqref{2.5} from Theorem \ref{th1}, together with H\"older inequality and the fact that $\min\big\{p,\frac{np}{n+p}\big\}>2$. The theorem has been established.
\ep

\smallskip

\section{Well-posedness of the direct Optical Tomography problem}   \label{Section3}

In this section we utilise the well-posedness results of Section \ref{Section2} to show that the direct problem of Fluorescent Optical Tomography is well posed.

\begin{theorem}[Well-posedness of the direct FOT problem] \label{th3} In the setting of Section \ref{Section1}, consider the boundary value problem \eqref{1.3} and suppose that the coefficients $\A,\K,\H,\xi, r, s, S, \ka,\la,\ga$ satisfy \eqref{1.4}-\eqref{1.6}, where $\Om \Subset \R^n$ is a domain with $C^1$ boundary and $n\geq 3$.

In addition, we also suppose that there exists $a_0>0$ such that
\beq 
\label{3.1}
\left\{
\begin{split}
&\A \in UC(\R^n;\R^{n\by n}_+),\ \ \, \K,\H \in L^{\infty}(\Om;\R^{2\by 2}),
\\
&  \ka \in C(\overline{\Om}),\ \ \ \ k_R,\ka \geq a_0 \, \text{ on }\, \Om, \ \ \ \ \ga,\la>0,
\end{split}
\right.
\eeq
where $``UC"$ is the space of bounded uniformly continuous functions. We further assume that for some $m>n$ we have
\beq
\label{3.2}
S \in L ^{\frac{nm}{2n+m}}(\Om;\R^{2}),\ \ \ s \in L^{\frac{m}{2}}(\p\Om;\R^2)
\eeq
and also
\[
p \, >\, \max\left\{n,\frac{2n}{n-2} \right\}.
\]
Then, for any $M>0$ and
\[
\xi \in C\big(\overline{\Om};[0,M]\big),
\]
the boundary value problem \eqref{1.3} has a unique weak solution
\[
(u,v) \, \in \, W^{1,\frac{m}{2}}(\Om;\R^{2}) \by W^{1,p}(\Om;\R^{2})
\]
which satisfies
\beq \label{3.3}
\left\{
\begin{split}
\phantom{\Bigg|}\!\! (a)\ \ & \int_\Om\Big[ \A_\xi :(\D u^\top \D \phi) \, +  \big(\K_\xi u-S\big)\cdot \phi \Big] \, \mathrm{d}\mL^n \,+ \int_{\p\Om}\big[ (\ga u-s)\cdot \phi\big] \, \mathrm{d}\mH^{n-1} =0 ,
\\
\!\! (b)\ \ & \int_\Om\Big[ \A_\xi :(\D v^\top \D \psi) \, +  \big(\K_\xi v- \xi \H u\big)\cdot \psi \Big] \, \mathrm{d}\mL^n \,+ \int_{\p\Om}\big[ \ga u\cdot \psi\big] \, \mathrm{d}\mH^{n-1} =0 ,
\end{split}
\right.
\eeq   
for all test functions
\[
(\phi,\psi) \, \in \, W^{1,\frac{m}{m-2}} (\Om;\R^{2}) \by W^{1,{\frac{p}{p-1}}} (\Om;\R^2). 
\]
In addition, there exists $C>0$ depending only on $M$, the coefficients, $p$ and the domain such that
\beq 
\label{3.4}
\left\{
\ \
\begin{split}
\!\!\!(a)\ \ \ \ \ \ \  \| u \|_{W^{1,\frac{m}{2}}(\Om)} +\, \| u \|_{L^m(\Om)} \, & \leq\, C \Big(\| S \|_{L ^{\frac{nm}{2n+m}}\!(\Om)} +\,\| s\|_{L^{\frac{m}{2}}\!(\p\Om)}\Big), \phantom{\Big|}
\\
\!\!\!(b)\ \ \ \ \, \| v \|_{W^{1,p}(\Om)}  +\, \| v \|_{C^{0,1-\frac{n}{p}}(\overline \Om)}\, & \leq\, C \, \| \xi \|_{L^\infty(\Om)} \| u \|_{L^m(\Om)} .
\end{split}
\right.
\eeq
\end{theorem}

\bp The goal is to apply Theorems \ref{th1}-\ref{th2}. To this end, we begin by showing that under \eqref{1.4}-\eqref{1.6} and \eqref{3.1} the diffusion matrix $\A_\xi$ satisfies the assumptions of these results. Since
\[
\A_\xi(x) \, =\, \A(x) \, +\, r(x,\xi)\mathrm I_n, \ \ \ r(x,t)\,=\, \frac{\la}{\ka(x)+t}
\]
and 
\[
\frac{\la}{\|\ka\|_{L^\infty(\Om)}+M} \, \leq\, r\big(x,\xi(x)\big) \, \leq\, \frac{\la}{a_0}
\]
the positive bounded uniformly continuous function 
\[
r\big(\cdot,\xi(\cdot)\big) : \ \ \overline{\Om} \ni x \mapsto r\big(x,\xi(x)\big) \in (0,\infty)
\]
can be extended to a positive bounded uniformly continuous function $\tilde r : \R^n \larrow (0,\infty)$ with the same upper and lower bounds as those of $r\big(\cdot,\xi(\cdot)\big)$. Then, since $\A$ is bounded and uniformly continuous on $\R^n$ with values in $\R^{n\by n}_+$, the matrix valued mapping 
\[
\tilde A \,:= \,A + \tilde r \mathrm I_n\ : \ \ \R^n \larrow \R^{n\by n}_+
\]
is bounded and uniformly continuous, valued in the positive matrices and its eigenvalues are uniformly bounded on $\overline{\Om}$ away from zero. Additionally, it is evident that $\K_\xi=\K+\xi \mathrm I_2 \in L^\infty(\Om;\R^{2\by 2})$ and that it satisfies the structural assumptions in \eqref{2.2}. Hence, by Theorems \ref{th1}-\ref{th2} applied to the Robin boundary value problem \eqref{1.3}(a)-\eqref{1.3}(c) for $p=m/2$ and noting that 
\[
\frac{m}{2}\, > \, \frac{n}{2} \,>\, \frac{2n}{n-2},
\]
for any $S \in L ^{\frac{nm}{2n+m}}(\Om;\R^{2})$ and any $s \in L^{\frac{m}{2}}(\p\Om;\R^2)$ there exists a unique solution $u \in W^{1,\frac{m}{2}}(\Om;\R^{2})$ satisfying \eqref{3.3}(a) for all $\phi \in W^{1,\frac{m}{m-2}} (\Om;\R^{2})$, as well as the estimate \eqref{3.4}(a). The only thing which is not already stated in the estimate \eqref{2.7} is the estimate on $\| u \|_{L^m(\Om)}$, which follows by the Sobolev inequalities.

Again by Theorems \ref{th1}-\ref{th2} applied to the Robin boundary value problem \eqref{1.3}(b)-\eqref{1.3}(d), there exists a unique solution $v \in W^{1,p}(\Om;\R^{2})$ satisfying \eqref{3.3}(b) for all $\psi \in W^{1,{\frac{p}{p-1}}} (\Om;\R^{2})$. Further, by applying \eqref{2.7}, by H\"older inequality we estimate
\[
\begin{split}
\| v \|_{W^{1,p}(\Om)} \, & \leq\, C \| \xi \H u\|_{L^{\frac{np}{n+p}}(\Om)}
\\
&\leq\, C  \| \xi \|_{L^{\infty}(\Om)}  \| u\|_{L^{\frac{np}{n+p}}(\Om)}
\\
& \leq\, C \| \xi \|_{L^\infty(\Om)} \| u\|_{L^{m}(\Om)},
\end{split}
\]
since $m> n$. The estimate \eqref{3.4}(b) therefore follows by the above estimate together with the Sobolev inequalities. The proof is complete. 
\ep

\smallskip
%\ms

\section{The inverse problem through PDE-constrained minimisation}   \label{Section4}

Now that the forward fluorescent optical tomography problem is understood, we proceed with the solvability of the inverse problem associated with \eqref{1.3}. Throughout this and subsequent sections we assume that the hypotheses of Theorem \ref{th3} are satisfied for a domain $\Om\Subset\R^n$ with $n\geq 3$ and which {\it from now is assumed to have $C^{1,1}$ regular boundary}.

Fix an integer $N\in \N$, $m >n$, $M,\de,\al>0$ and $p >\max\left\{n,\frac{2n}{n-2} \right\}$. Consider Borel sets
\beq \label{4.1}
\big\{\Ga_1,...,\Ga_N\big\}  \sub \p\Om
\eeq
and light sources 
\beq  \label{4.2}
\big\{S_1,...,S_N \big\} \sub L ^{\frac{nm}{2n+m}}(\Om;\R^{2})\ , \ \ \ \ \big\{s_1,...,s_N\big\} \sub L^{\frac{m}{2}}(\p\Om;\R^{2})
\eeq
in the interior and on the boundary respectively. Let also 
\beq
 \label{4.3}
\big\{v^\de_1,...,v^\de_N\big\}  \sub L^{\infty}(\p\Om;\R^{2})
\eeq
be predicted approximate values of the solution $v$ of \eqref{1.3}(b)-\eqref{1.3}(d) on the boundary $\p\Om$, at noise (error) level $\de$. Suppose that for any $i\in\{1,...,N\}$, the pair $(u_i,v_i)$ solves \eqref{1.3} with data $(S_i,s_i,\xi)$. For the $N$-tuple of solutions $(u_1,...,u_N;v_1,...,v_N)$, we will be using the notation
\[
\big(\vec u,\vec v \, \big) \, \in \, W^{1,\frac{m}{2}}(\Om;\R^{2\by N}) \by W^{1,p}(\Om;\R^{2 \by N})
\]
and understand the vectors of solutions $(u_i)_{i=1...N}$ and $(u_i)_{i=1...N}$ as being $\R^{2\by N}$-valued. Similarly, we will understand the corresponding vectors of test functions as
\[
\big(\vec \phi,\vec \psi \, \big) \, \in \, W^{1,\frac{m}{m-2}}(\Om;\R^{2\by N}) \by W^{1,{\frac{p}{p-1}}}(\Om;\R^{2 \by N}).
\]
The goal of the inverse problem associated with \eqref{1.3} is to {\it determine a non-negative $\xi \in L^p(\Om,[0,\infty))$ such that the errors $\big|(v_i - v^\de_i)\big|_{\Ga_i}\big|$ which describe the misfit between the predicted approximate solution and the actual solution become as small as possible}. We will minimise the error in $L^\infty$ by means of approximations in $L^p$ for large $p$ and then take the limit $p\to \infty$. The benefit of minimisation in $L^\infty$ is that one can achieve uniformly small error rather than on average. Since no reasonable error functional is coercive in the admissible class of $N$-tuples of PDE solutions without additional constraints, we add an extra Tykhonov-type regularisation term $ \al \| \xi \|$ for a small parameter $\al>0$ and some appropriate norm. 

In view of the above observations, we define for $p>\max\left\{n,\frac{2n}{n-2} \right\}$ the functional
\beq 
\label{4.4}
\E_p\big(\vec u,\vec v ,\xi \big) \,:=\, \sum_{i=1}^N \big\|v_i - v^\de_i\big\|_{\dot L^p(\Ga_i)}  +\, \al \big\| \D^2 \xi \big\|_{\dot L^m(\Om)}, \ \ \ \ (\vec u,\vec v \,,\xi \big) \in \mathfrak{X}^{p}(\Om)
\eeq
and its limiting counterpart
\beq
\label{4.5}
\E_\infty\big(\vec u,\vec v ,\xi \big) \,:=\, \sum_{i=1}^N \big\|v_i - v^\de_i\big\|_{L^\infty(\Ga_i)}  +\, \al \big\| \D^2 \xi \big\|_{\dot L^m(\Om)} \ \ \ \ (\vec u,\vec v \,,\xi \big) \in \mathfrak{X}^{\infty}(\Om),
\eeq
where the dotted $\dot L^p$-functionals are the next regularisations of the respective norms
\beq
\label{4.6}
 \|g \|_{\dot L^p(\Ga_i)} := \left(\, \av_{\Ga_i}(|g|_{(p)})^p\, \mathrm d \mH^{n-1} \!\right)^{\!\!1/p}\ , \ \ \ \ 
\|f \|_{\dot L^m(\Om)} := \left(\, \av_{\Om}(|f|_{(m)})^m\, \mathrm d \mL^{n}\! \right)^{\!\!1/m}
\eeq
and $|\cdot |_{(p)}$ is a regularisation of the Euclidean norm away from zero in the corresponding space, given by
\beq
\label{4.7}
| \cdot |_{(p)} \,:=\, \sqrt{| \cdot |^2 + p^{-2}}.
\eeq
 {The slashed integral symbolises the average with respect to either the Lebesgue measure $\mL^n$ or the Hausdorff measure $\mH^{n-1}$.} The respective admissible classes $\mX^p(\Om)$ and $\mX^\infty(\Om)$ are defined by
\beq
\label{4.8}
\mX^p(\Om)\,:=\, \left\{ 
\begin{array}{l}
 (\vec u,\vec v \,,\xi \big) \in \mW^{p}(\Om) :  \text{ for any } i \in\{1,...,N\},\ (u_i,v_i,\xi) \text{ satisfies }\ms\ms
\\
\hspace{80pt} 0\leq \xi \leq M \, \text{ a.e. on $\Om$}
\\
\text{and}\ms
\\
\ \  \left\{\ \ 
\begin{array}{lll}
(a)_i\ \ & \ \  -\div(\A_\xi  \hspace{1pt} \overset{_\bullet}{\phantom{.}}  \D u_i)\,+\, \K_\xi u_i\, =\, S_i, & \ \ \text{ in }\Om, \ms
\\
(b)_i&\ \     -\div(\A_\xi \hspace{1pt} \overset{_\bullet}{\phantom{.}}  \D v_i)\,+\, \K_\xi v_i\, =\, \xi \H u_i, & \ \ \text{ in }\Om,  \ms
\\
(c)_i&  \ \ \ \ (\A_\xi  \hspace{1pt} \overset{_\bullet}{\phantom{.}}  \D u_i)\hspace{1pt} \mathrm n\,+\, \ga u_i\, =\, s_i, & \ \ \text{ on }\p\Om,  \ms
\\
(d)_i& \ \ \ \  (\A_\xi \hspace{1pt} \overset{_\bullet}{\phantom{.}}  \D v_i)\hspace{1pt} \mathrm n\,+\, \ga v_i\, =\, 0, & \ \ \text{ on }\p\Om,
\end{array}
\right.
\\
\\
 \text{for }\A,\K,\H,S_i,s_i,\ka,\xi,\la,\ga,p\text{ satisfying the hypotheses \eqref{1.4}}
\\
\text{\eqref{1.5}, \eqref{1.6}, \eqref{3.1} and \eqref{4.1}-\eqref{4.3}}
\end{array}\!\!
\right\}
\eeq
and
\beq
\label{4.9}
\mX^\infty(\Om)\, :=\, \bigcap_{n<p<\infty} \mX^p(\Om),
\eeq
whilst the  {Banach space $\mW^{p}(\Om)$ involved in the definition of the admissible class $\mX^p(\Om)$ is}
\beq
\label{4.10}
 {\phantom{\Big|}\mW^p(\Om)\,:=\, W^{1,\frac{m}{2}}(\Om;\R^{2\by N}) \by W^{1,p}(\Om;\R^{2 \by N}) \by W^{2,m}(\Om).}
\eeq
 {Note that $\mX^\infty(\Om)$ is a subset of a Frech\'et space, rather than of a Banach space, but this will not cause any added difficulties. In particular, note that each $\mX^p(\Om)$ is subset of the Banach space $\mW^p(\Om)$, but $\mX^\infty(\Om)$ is {\it not} ``$\mX^p(\Om)$ for $p=\infty$" because our a priori estimates are not true in $W^{1,\infty}(\Om)$, but only in $W^{1,p}(\Om)$ for any $p \in[2,\infty)$ with a constant $C=C(p)$ depending on $p$. (This reminiscent to the classical strict inclusion $L^\infty(\Om) \subsetneqq \bigcap_{1<p<\infty}L^p(\Om)$ which is true when $\mL^n(\Om)<\infty$: any $f\in \bigcap_{1<p<\infty}L^p(\Om)$ satisfies $\|f\|_{L^p(\Om)}\leq C(p)<\infty$, but there exist $f$'s with $C(p)\to \infty$ as $p\to \infty$ for which $\|f\|_{L^\infty(\Om)}=\infty$.)}

\begin{remark}It might be quite surprising that in the Tikhonov term we include the $L^m$ norm of the Hessian of $\xi$, rather than as one would expect the $L^m$ norm of either the gradient or $\xi$ itself. It turns out that one cannot regularise enough by adding ``$+\al \| \xi \|_{\dot L^m(\Om)}$" to obtain minimisers (this would be redundant anyway because of the unilateral constraint). On the other hand, by adding ``$+\al \| \D\xi \|_{\dot L^m(\Om)}$", one can indeed recover all the results up to and including Section \ref{Section5}, but not the results of Section \ref{Section6}, as we cannot obtain the variational inequalities in $L^\infty$ without higher regularity in the coefficients of the PDE systems in \eqref{4.8} due to the emergence of  quadratic gradient terms.
\end{remark}

The main result in this section concerns the existence of $\E_p$-minimisers in the admissible class $\mX^p(\Om)$, the existence of $\E_\infty$-minimisers in the admissible class $\mX^\infty(\Om)$ and the approximability of the latter by the former as $p\to \infty$.

\begin{theorem}[$\E_\infty$-error minimisers, $\E_p$-error minimisers and convergence as $p\to\infty$] \label{th4}

\ms

\noi {\rm (A)} The functional $\E_p$  given by \eqref{4.4}  {has a unique constrained minimiser} $(\vec u_p,\vec v_p ,\xi_p)$ in the admissible class $\mX^p(\Om)$:
\beq \label{4.15A}
\E_p\big(\vec u_p,\vec v_p ,\xi_p \big)\, =\, \inf\Big\{\E_p\big(\vec u,\vec v ,\xi \big)\, : \ \big(\vec u,\vec v ,\xi \big) \in \mX^p(\Om) \Big\}.
\eeq
{\rm (B)} The functional $\E_\infty$ given by \eqref{4.5} has a constrained minimiser $(\vec u_\infty,\vec v_\infty ,\xi_\infty)$ in the admissible class $\mX^\infty(\Om)$:
\beq  \label{4.15B}
\E_\infty\big(\vec u_\infty,\vec v_\infty ,\xi_\infty \big)\, =\, \inf\Big\{\E_\infty\big(\vec u,\vec v ,\xi \big)\, : \ \big(\vec u,\vec v ,\xi \big) \in \mX^\infty(\Om) \Big\}.
\eeq
Additionally, there exists a subsequence of indices $(p_j)_1^\infty$ such that the sequence of respective $\E_{p_j}$-minimisers $\big(\vec u_{p_j},\vec v_{p_j} ,\xi_{p_j} \big)$ satisfies
\beq
\label{4.16}
\left\{ \ \
\begin{array}{ll}
\vec u_{p} \weak \vec u_\infty, & \text{ in }W^{1,\frac{m}{2}}(\Om;\R^{2\by N}), \ms
\\
\vec u_{p} \larrow \vec u_\infty, & \text{ in }L^{\frac{m}{2}}(\Om;\R^{2\by N}), \ms
\\
\vec v_{p} \weak \vec v_\infty, & \text{ in }W^{1,q}(\Om;\R^{2\by N}), \text{ for all }q\in (1,\infty),\ms
\\
\vec v_p \larrow \vec v_\infty, & \text{ in }C(\overline{\Om};\R^{2\by N}),  \ms
\\
\xi_p \weak \xi_\infty, & \text{ in }W^{2,m}(\Om),
\ms
\\
\xi_p \larrow \xi_\infty, & \text{ in }C^1(\overline{\Om}),
\end{array}
\right.
\eeq
as $p_j\to \infty$. Further
\beq \label{4.14a}
\E_\infty\big(\vec u_\infty,\vec v_\infty ,\xi_\infty\big)\, =\, \lim_{p_j\to\infty}\E_p\big(\vec u_{p},\vec v_{p} ,\xi_{p} \big).
\eeq
\end{theorem}

The proof of Theorem \ref{th4} is a consequence of the next two propositions, utilising the direct method of Calculus of Variations (\cite{D}).

\begin{proposition}[$\E_p$-minimisers] \label{pr4} In the setting of Theorem \ref{th4}, the functional $\E_p$  {has unique constrained minimiser} $\big(\vec u_p,\vec v_p ,\xi_p \big)$ in the admissible class $\mX^p(\Om)$, as claimed in \eqref{4.15A}.
\end{proposition}

\bp Let us begin by noting that $\mX^p(\Om) \neq \emptyset$, and, as we will show right next, in fact it is a weakly closed subset of the reflexive Banach space $\mW^p(\Om)$ with cardinality greater or equal to that of $L^p(\Om)$. Next note that there is an a priori energy bound for the infimum of $\E_p$, in fact uniform in $p$. Indeed, for each $i\in\{1,...,N\}$ let $(u_{0i},v_{0i})$ be the solution to \eqref{1.3} with $\xi \equiv 0$ and sources $(S_i,s_i)$ as in \eqref{4.2}. Then, by Theorem \ref{th3}, we have $v_{0i} \equiv 0$. Therefore, by \eqref{4.5}-\eqref{4.6} we infer that
\[
\big(\vec u_0,\vec 0 ,0 \big) \, \in \, \mX^p(\Om)
\]
for all $p \in [n,\infty]$, and also, by H\"older inequality and \eqref{4.3}-\eqref{4.7}, we obtain
\[
\begin{split}
\E_p\big(\vec u_0,\vec 0 ,0 \big) \, & \leq\, \frac{N+1}{p} \,+\, E_\infty\big(\vec u_0,\vec 0 ,0 \big)
\\
& \leq\,  \frac{N+1}{n} \,+\, \sum_{i=1}^N \big\| v^\de_i\big\|_{L^\infty(\Ga_i)}
\\
&<\, \infty.
\end{split}
\]
Consider now (for fixed $p$) a minimising sequence $(\vec u^j,\vec v^{\, j},\xi^{\, j}\,)_{j=1}^\infty$ of $\E_p$ in $\mX^p(\Om)$. Then, for all large enough $j\in\N$ we have 
\[
0\,\leq\, \E_p\big(\vec u^{\,j} ,\vec v^{\, j}, \xi^{\, j}\big) \, \leq\, 1\,+ \, \frac{N+1}{n} \,+\, \sum_{i=1}^N \big\| v^\de_i \big\|_{L^\infty(\Ga_i)}.
\]
By Theorem \ref{th3}, we have the estimates
\beq
\label{4.13}
\left\{
\ \
\begin{split}
\phantom{_\big|}  \big\| \vec u^{\,j} \big\|_{W^{1,\frac{m}{2}}(\Om)}   \, & \leq\, C \, \sum_{i=1}^N\Big(\| S_i \|_{L ^{\frac{nm}{2n+m}}\!(\Om)} +\,\| s_i\|_{L^{\frac{m}{2}}\!(\p\Om)}\Big), \phantom{\Big|}
\\
\big\| \vec v^{\, j}\big\|_{W^{1,p}(\Om)} \, & \leq\, C \, M \, \big\| \vec u^{\, j}\big\|_{W^{1,\frac{m}{2}}(\Om)} .
\end{split}
\right.
\eeq
By the above and \eqref{4.4}, we have the estimate
\[
\big\| \D^2 \xi^{\,j} \big\|_{L^m(\Om)} \, \leq\, \big\| \D^2 \xi^{\,j} \big\|_{\dot L^m(\Om)} \, \leq\, \frac{1}{\al} \bigg( 1 + \frac{N+1}{n}+ \sum_{i=1}^N \big\| v^\de_i\big\|_{L^\infty(\Ga_i)}\bigg).
\]
Further, in view of the unilateral constraint, we readily have
\[
\| \xi^{\, j} \|_{L^\infty(\Om)} \, \leq\, M .
\]
By the $C^{1,1}$ regularity of $\p\Om$, the interpolation inequalities in the Sobolev space $W^{2,m}(\Om)$ (see e.g.\ \cite[Theorem 7.28, p.173]{GT}) imply the existence of $C>0$ independent of $j$ such that
\[
\big\| \D \xi^{\,j} \big\|_{L^m(\Om)}  \leq\, C\big\| \D^2 \xi^{\,j} \big\|_{L^m(\Om)}  +\, \| \xi^{\, j} \|_{L^m(\Om)}.
\]
Thus, the above estimates yield the uniform bound
\beq
\label{4.14}
\sup_{j\in\N}\, \| \xi^{\, j}\|_{W^{2,m}(\Om)} \, <\, \infty .
\eeq
By the estimates \eqref{4.13}-\eqref{4.14} and standard weak and strong compactness arguments, there exists a weak limit
\[
(\vec u_p,\vec v_p , \xi_p \big) \, \in \mW^p(\Om)
\]
and a subsequence $(j_k)_1^\infty$ such that along which we have
\[
\left\{ \ \
\begin{array}{ll}
\vec u^{\, j} \weak \vec u_p, & \text{ in }W^{1,\frac{m}{2}}(\Om;\R^{2\by N}), \ms
\\
\vec u^{\, j} \larrow \vec u_p, & \text{ in }L^{\frac{m}{2}}(\Om;\R^{2\by N}), \ms
\\
\vec v^{\, j} \weak \vec v_p, & \text{ in }W^{1,p}(\Om;\R^{2\by N}), \ms
\\
\vec v^{\, j} \larrow \vec v_p, & \text{ in }C(\overline{\Om};\R^{2\by N}),  \ms
\\
\xi^{\, j} \weak \xi_p, & \text{ in }W^{2,m}(\Om), \ms
\\
\xi^{\, j} \larrow \xi_p, & \text{ in }C^1(\overline{\Om}),
\end{array}
\right.
\]
as $j_k \to \infty$. We note that in this paper we will utilise the common practice of passing to subsequences as needed without perhaps explicit relabelling of the new subsequences. To show that in fact the limit $(\vec u_p,\vec v_p , \xi_p \big)$  lies in the admissible constrained class $\mX^p(\Om)$, we argue as follows. Note now that the pointwise constraint
\[
0 \, \leq\, \xi^{\,j}   \leq\, M \ \text{ a.e. on }\Om,
\]
is weakly closed in $W^{2,m}(\Om)$, namely the set
\beq
\label{4.15}
W^{2,m}(\Om;[0,M])\,=\, \Big\{\eta \in L^p(\Om):\, \ 0 \, \leq\, \eta\,   \leq\, M \ \text{ a.e. on }\Om\Big\}
\eeq
is weakly closed. This is an immediate consequence of the strong compactness of the embedding of $W^{2,m}(\Om)$ into $C^1(\overline{\Om})$. We thus infer that 
\[
(\vec u_p,\vec v_p , \xi_p \big) \in \mX^p(\Om)
\]
by passing to the weak limit in the equations $(a)_i-(d)_i$ defining \eqref{4.8}, which is possible due to the modes of convergence the minimising sequence satisfies.

We finally show that the weak limit $(\vec u_p,\vec v_p , \xi_p \big) \in \mX^p(\Om)$ is indeed the minimiser of $\E_p$ over the same space. To this end, note that for any $\al>0$ the nonlinear functional $\al \| \D^2(\, \cdot\, ) \|_{\dot L^m(\Om)}$ is convex and strongly continuous on the reflexive space $W^{2,m}(\Om)$, by \eqref{4.6}-\eqref{4.7}. Therefore, it is weakly lower semi-continuous on the same space. Similarly, for each $i\in\{1,...,N\}$ the functional $\big\|\cdot - v^\de_i\big\|_{\dot L^p(\Ga_i)}$ is strongly continuous on $L^p(\Ga_i)$. Hence, we conclude that
\[
\begin{split}
\E_p\big(\vec u_p,\vec v_p ,\xi_p \big) \,&=\, \sum_{i=1}^N \big\|v_{pi} - v^\de_i\big\|_{\dot L^p(\Ga_i)}  +\, \al \big\| \D^2 \xi_p \big\|_{\dot L^m(\Om)} 
\\
&\leq \liminf_{j_k\to \infty} \left\{ \sum_{i=1}^N \big\|v^{\, j}_{i} - v^\de_i\big\|_{\dot L^p(\Ga_i)}  +\, \al \big\| \D^2 \xi^{\, j} \big\|_{\dot L^m(\Om)} \right\}
\\
&= \liminf_{j_k\to \infty} \E_p\big(\vec u^{\,j} ,\vec v^{\, j}, \xi^{\, j}\big)
\\
&=\, \inf\big\{\E_p :\ \mX^p(\Om) \big\}.
\end{split}
\]
 {Finally, note that for any fixed $p$, the $p$-minimiser $(\vec u_{p},\vec v_{p} ,\xi_{p} \big)$ is unique. On the one hand, the functional $\E_p$ is strictly convex in $\xi$ (and also in $\vec v$, but constant in $\vec u$). On the other hand, for any $\xi$, there exists a unique pair $(\vec u,\vec v\big)$ satisfying the PDE constraints, as a consequence of Theorem \ref{th3}.} The proposition ensues.
\ep

Our next result below concerns the existence of minimisers for the $L^\infty$-error functional and approximation of those by minimisers of $L^p$ functionals, completing the proof of Theorem \ref{th4}.

\begin{proposition}[$\E_\infty$-minimisers] \label{pr5} In the above setting, the functional $\E_\infty$ given by \eqref{4.5} has a constrained minimiser $(\vec u_\infty,\vec v_\infty ,\xi_\infty)$ in the admissible class $\mX^\infty(\Om)$, as claimed in \eqref{4.15B}. 

Additionally, there exists a subsequence of indices $(p_j)_1^\infty$ such that the sequence of respective $\E_{p_j}$-minimisers $\big(\vec u_{p_j},\vec v_{p_j} ,\xi_{p_j} \big)$ constructed in Proposition \ref{pr4} satisfy \eqref{4.16} as $p_j\to \infty$. Further the energies converge as claimed in \eqref{4.14a}.
\end{proposition}

\bp We essentially continue from the proof of Proposition \ref{pr4}. The energy bound $(\vec u_0,\vec 0,0)$ constructed therein is uniform in $p$ and also, in view of \eqref{4.8}-\eqref{4.9} we have $(\vec u_0,\vec 0,0) \in\mX^\infty(\Om)$. Fix now $q>n$ and consider large enough $p\geq q$. Then, by H\"older inequality and minimality, we have the bound
\[
\begin{split}
\E_q\big(\vec u_{p},\vec v_{p} ,\xi_{p} \big)\, &\leq\, \E_p\big(\vec u_{p},\vec v_{p} ,\xi_{p} \big)
\\
 & \leq\, \E_p(\vec u_0,\vec 0,0)
 \\
 & \leq\, \frac{N+1}{n} + \, \E_\infty(\vec u_0,\vec 0,0)
 \\
  & \leq \, \frac{N+1}{n}+ \sum_{i=1}^N \big\| v^\de_i\big\|_{L^\infty(\Ga_i)},
\end{split}
\]
which is uniform in $p$. By the above estimate, we have
\[
\big\| \D^2 \xi_p \big\|_{L^m(\Om)} \, \leq\, \big\| \D^2\xi_p \big\|_{\dot L^m(\Om)} \, \leq\, \frac{1}{\al} \bigg(\frac{N+1}{n}+ \sum_{i=1}^N \big\| v^\de_i\big\|_{L^\infty(\Ga_i)} \bigg).
\]
On the other hand, by the unilateral pointwise constraint, we immediately have
\[
0 \, \leq\,  \xi_p  \leq\, M \  \text{ on }\overline{\Om} .
\]
Hence, by the interpolation inequalities in $W^{2,m}(\Om)$, we deduce the uniform bound
\beq
\label{4.14}
\sup_{p\geq n}\, \| \xi_p\|_{W^{2,m}(\Om)} \, <\, \infty .
\eeq
By the above estimates, by Theorem \ref{th3} (see \eqref{4.13}) and by standard weak and strong compactness arguments together with a diagonal argument, there exists a limit
\[
(\vec u_\infty,\vec v_\infty , \xi_\infty \big) \, \in \bigcap_{n<q<\infty}\mW^q(\Om)
\]
and a subsequence $(p_j)_1^\infty$ such that the modes of convergence in \eqref{4.16} hold true as $p_j \to \infty$. Further, by passing to the limit as $p_j\to \infty$ in the equations $(a)_i-(d)_i$ forming the admissible class \eqref{4.8} and the closed unilateral pointwise constraint $0\leq \xi_p \leq M$, we see that in fact the limit $(\vec u_\infty,\vec v_\infty , \xi_\infty \big)$ lies in the admissible class $\mX^\infty(\Om)$. It remains to show that $(\vec u_\infty,\vec v_\infty , \xi_\infty \big)$ is indeed a minimiser of $\E_\infty$ and that the energies converge as claimed. To this end, fix an arbitrary $(\vec u,\vec v, \xi \big) \in \mX^\infty(\Om)$. Since $p_j\geq q$ for any $q>1$ and large enough $j\in\N$, by minimality we have
\[
\begin{split}
\E_\infty\big(\vec u_\infty,\vec v_\infty ,\xi_\infty \big) \,&=\, \lim_{q\to \infty} \E_q\big(\vec u_\infty,\vec v_\infty ,\xi_\infty \big)
\\
&\leq\, \liminf_{q\to \infty} \Big(\liminf_{p_j\to \infty}\, \E_q\big(\vec u_p,\vec v_p ,\xi_p \big)  \Big)
\\
&\leq\, \liminf_{p_j\to \infty} \E_p\big(\vec u_p,\vec v_p ,\xi_p \big) 
\\
&\leq\, \limsup_{p_j\to \infty} \E_p\big(\vec u_p,\vec v_p ,\xi_p \big)  
\\
&\leq\, \lim_{p_j\to \infty} \E_p\big(\vec u,\vec v,\xi \big)  
\\
&=\,  \E_\infty\big(\vec u,\vec v ,\xi \big) , 
\end{split}
\]
for any $(\vec u,\vec v, \xi \big) \in \mX^\infty(\Om)$. Hence $\big(\vec u_\infty,\vec v_\infty ,\xi_\infty \big)$ is a minimiser of $\E_\infty$ over $\mX^\infty(\Om)$. The particular choice of $(\vec u,\vec v, \xi \big):=\big(\vec u_\infty,\vec v_\infty ,\xi_\infty \big)$ in the above inequality yields
\[
\lim_{p_j\to \infty} \E_p\big(\vec u_p,\vec v_p ,\xi_p \big) \, =\, \E_\infty\big(\vec u_\infty,\vec v_\infty ,\xi_\infty \big). 
\]
The proof of Proposition \ref{pr5} is now complete.
\ep

\ms

\section{Kuhn-Tucker theory and Lagrange multipliers for the $p$-error}  \label{Section5}

In this section we return the $L^p$-minimisation problem \eqref{4.15A} solved in Theorem \ref{th4} for finite $p<\infty$ (Section \ref{Section4}). Given the presence of both PDE and unilateral constraints, in general one cannot have an Euler-Lagrange equation, but an one-sided variational inequality with Lagrange multipliers. The goal here is to derive the relevant variational inequality associated with \eqref{4.15A}. The main result is therefore the following.

\begin{theorem}[The variational inequalities in $L^p$] \label{th7}

In the setting of Section \ref{Section4} and under the same assumptions, for any $p>\max\{n,2n/(n-2)\}$, there exist Lagrange multipliers
\[
\big(\vec \phi_p, \vec \psi_p \big) \ \in W^{1,\frac{m}{m-2}} (\Om;\R^{2\by N}) \by W^{1,\frac{p}{p-1}} (\Om;\R^{2\by N})
\]
associated with the constrained minimisation problem \eqref{4.15A} for $\E_p$ in the admissible class $\mX^p(\Om)$, such that the constrained minimiser $\big(\vec u_p, \vec v_p, \xi_p\big) \in \mX^p(\Om)$ satisfies the next three relations:
\beq
\label{5.20}
\begin{split}
& \frac{\al m}{p}\int_\Om (\D^2\eta-\D^2\xi_p) : \mu(\D^2\xi_p)\, \mathrm{d}\mL^n 
 \, \geq
\\
& \sum_{i=1}^N  \int_\Om (\eta- \xi_p) \bigg\{ \!\! -\big(\H  u_{pi}  \big)  \cdot \psi_{pi}\,+\, \dot r(\cdot, \xi_p)\Big[  \D u_{pi} : \D \phi_{pi} 
\\
& \hspace{70pt} + \D v_{pi} :\D \psi_{pi} \Big] +  u_{pi} \cdot \phi_{pi} + v_{pi} \cdot \psi_{pi}    \bigg\} \, \mathrm{d}\mL^n
\end{split}
\eeq
for any $\eta \in W^{2,m}(\Om,[0,M]) $; further,
\beq
\label{5.21}
\begin{split}
 \int_{\p\Om} \vec w : \mathrm d [\vec \nu_p(\vec v_p)] \, =\,\sum_{i=1}^N & \bigg\{  \int_\Om \Big[ \A_{\xi_p} :(\D w_i^\top \D \psi_{pi}) \, +  \big(\K_{\xi_p} w_i\big) \cdot \psi_{pi} \Big] \, \mathrm{d}\mL^n 
\\
& \ + \int_{\p\Om} (\ga w_i)\cdot \psi_{pi} \, \mathrm{d}\mH^{n-1} \bigg\},
\end{split}
\eeq
for any $\vec w  \in W^{1,p} (\Om;\R^{2\by N})$, and finally
\beq
\label{5.22}
\begin{split}
  \sum_{i=1}^N & \bigg\{\int_\Om\Big[ \A_{\xi_p} :(\D z_i^\top \D \phi_{pi}) \, + (\K_{\xi_p} z_i)\cdot \phi_{pi} \Big] \, \mathrm{d}\mL^n \,+ \int_{\p\Om}(\ga z_i)\cdot \phi_{pi} \, \mathrm{d}\mH^{n-1} \bigg\}
\\
& \ \ =\ \sum_{i=1}^N \int_\Om   \xi_p \big(\H z_i \big) \cdot \psi_{pi} \, \mathrm{d}\mL^n  , 
\end{split}
\eeq
for any $ \vec z  \in W^{1,\frac{m}{2}} (\Om;\R^{2\by N})$. 

\ms

In the relations \eqref{5.20}-\eqref{5.21}, $\mu(V) $ is defined for any $V \in L^m(\Om;\R^{n\by n})$ as
\beq 
\label{5.3}
\mu (V)\, :=\, \frac{ (|V|_{(m)})^{m-2} \, V}{ \mL^{n}(\Om) \big(\|V \|_{\dot L^m(\Om)}\big)^{m-1}}  ,
\eeq
$\dot r(x,t)$ symbolises the partial derivative of $r$ with respect to $t$ (recall \eqref{1.4}) and $\vec \nu_p(\vec v)$ is a $\vec v$-dependent $\R^{2\by N}$-valued matrix Radon measure in $\mM(\p\Om;\R^{2\by N})$ given by
\beq 
\label{5.2}
\vec \nu_p(\vec v)\, :=\, \sum_{i=1}^N \Bigg(\frac{\big(|v_i-v^\de_i|_{(p)}\big)^{p-2} (v_i-v^\de_i)}{ \mH^{n-1}(\Ga_i) \big(\big\|v_i-v^\de_i\big\|_{\dot L^p(\Ga_i)}\big)^{p-1}}  \ot e_i \Bigg)  \mH^{n-1}\LL_{\Ga_i} .
\eeq
\end{theorem}

Note that $\vec \nu_p(\vec v)$ is absolutely continuous with respect to the Hausdorff measure $\mH^{n-1}\LL_{\p\Om}$. The notation $\{e_1,...,e_N\}$ symbolises the standard Euclidean basis of $\R^N$ and ``$:$" symbolises the standard inner product in $\R^{2\by N}$. Additionally, one may trivially compute that 
\[
\dot r(x,t) \, =\, -\frac{\la}{(\ka(x)+t)^2}.
\]
\begin{remark} The reason that we obtain the three different relations \eqref{5.20}-\eqref{5.22} of which one is inequality and two are equations can be explained as follows. If one ignores the PDE constraints defining  \eqref{4.8} (which give rise to the Lagrange multipliers), then the admissible class is in fact the Cartesian product of three sets, two of which are vector spaces (spaces for $\vec u$ and $\vec v$), and one is a convex set (space of $\xi$), see \eqref{5.10} that follows. Hence, since the unilateral constraint is only for $\xi$, the variational inequality is only for this variable. The decoupling of the relations is merely a consequence of linear independence.
\end{remark}

The proof of Theorem \ref{th7} consists of several particular sub-results. We begin with computing the Frech\'et derivative of the functional $\E_p$.

\begin{lemma} \label{pr7} The functional $\E_p : \mW^p(\Om) \larrow \R$ given by \eqref{4.4}-\eqref{4.10} is Frech\'et differentiable and its derivative
\[
\mathrm d \hspace{0.5pt} \E_p \ : \ \ \mW^p(\Om) \larrow \big(\mW^p(\Om)\big)^*\ , \ \ \ \ (\vec u, \vec v,\xi) \mapsto \big(\mathrm d \hspace{0.5pt} \E_p\big)_{ (\vec u, \vec v,\xi)},
\]
is given by 
\beq 
\label{5.1}
\big(\mathrm d \hspace{0.5pt} \E_p\big)_{ (\vec u, \vec v,\xi)} (\vec z, \vec w,\eta)  \,=\, p\int_{\p\Om} \vec w : \mathrm d [\vec \nu_p(\vec v)]\, + \, \al m \int_\Om  \D^2\eta : \mu(\D^2 \xi)\, \mathrm d \mL^n
\eeq
for any $(\vec u, \vec v,\xi), (\vec z, \vec w,\eta) \in \mW^p(\Om)$. 
\end{lemma}

\bp The Frech\'et differentiability of $\E_p$ follows from standard results on the geometry of Banach spaces and the $p$-regularisations of the norms, given by \eqref{4.6}-\eqref{4.7}. To compute the Frech\'et derivative, we use Gateaux differentiation. To this end, fix $(\vec u, \vec v,\xi), (\vec z, \vec w,\eta) \in \mW^p(\Om)$. Then, we have
\[
\begin{split}
 \big(\mathrm d \hspace{0.5pt} \E_p\big)_{ (\vec u, \vec v,\xi)} (\vec z, \vec w,\eta)  \, & =\, \frac{\mathrm d}{\mathrm d \e}\Big|_{\e=0}\E_p \Big((\vec u, \vec v,\xi)+\e (\vec z, \vec w,\eta)\Big)
\\
&=\, \sum_{i=1}^N \frac{\mathrm d}{\mathrm d \e}\Big|_{\e=0}\left( \, \av_{\Ga_i}\big(\big| v_i+\e w_i -v^\de_i \big|_{(p)} \big)^p\, \mathrm d \mH^{n-1} \! \right)^{\! \frac{1}{p}}
\\
& \ \ \ \ + \, \al \, \frac{\mathrm d}{\mathrm d \e}\Big|_{\e=0}\left( \, \av_{\Om}\left(\big|\D^2\xi+\e \D^2\eta\big|_{(m)}\right)^m\, \mathrm d \mL^{n} \! \right)^{\! \frac{1}{m}}
\end{split}
\]
which by the chain rule yields
\[
\begin{split}
 \big(\mathrm d \hspace{0.5pt} \E_p\big)_{ (\vec u, \vec v,\xi)} (\vec z, \vec w,\eta)  \, 
& = \, p\, \sum_{i=1}^N  \left( \, \av_{\Ga_i}\big(\big| v_i -v^\de_i \big|_{(p)} \big)^p\, \mathrm d \mH^{n-1} \! \right)^{\! \frac{1}{p}-1} \centerdot
\\
&\ \ \ \ \,  \centerdot \av_{\Ga_i}\big(|v_i-v^\de_i|_{(p)}\big)^{p-2} (v_i-v^\de_i) \cdot w_i\, \mathrm d \mH^{n-1} 
\\
  + \, \al m & \left( \, \av_{\Om}\big(|\D^2\xi|_{(m)}\big)^m\, \mathrm d \mL^{n} \! \right)^{\! \frac{1}{m} -1}
{\av_{\Om}}\big(\big|\D^2\xi\big|_{(m)}\big)^{m-2}\, \D^2\xi : \D^2\eta \, \mathrm d \mL^{n}.
\end{split}
\]
Hence, \eqref{5.1} follows in view of the definitions \eqref{5.3}-\eqref{5.2}. The lemma ensues.
\ep

In order to derive the appropriate variational inequality that any minimiser as in \eqref{4.15A} satisfies, we need to define a map which incorporates the PDE constraints of the admissible class $\mX^p(\Om)$ in \eqref{4.8}.
We define
\beq \label{5.6}
\G \ : \ \ \mW^p(\Om) \larrow \Big[\big( W^{1,\frac{m}{m-2}} (\Om;\R^{2})\big)^*  \by \big(W^{1,\frac{p}{p-1}} (\Om;\R^{2})\big)^*\Big]^N
\eeq
by setting
\beq  \label{5.7}
\begin{split}
& \left\langle \G(\vec u, \vec v, \xi )  ,\, (\vec \phi, \vec \psi\, )\right\rangle\, := \, \left[
\begin{array}{c}
 \Big\langle \G^1_1(\vec u, \vec v, \xi )  ,\,  \phi_1 \Big\rangle \ms
 \\
  \Big\langle \G^2_1(\vec u, \vec v, \xi )  ,\,  \psi_1 \Big\rangle 
 \\
 \vdots
 \\
 \Big\langle \G^1_N(\vec u, \vec v, \xi )  ,\, \phi_N \Big\rangle \ms
 \\
 \Big\langle \G^2_N(\vec u, \vec v, \xi )  ,\,  \psi_N\Big\rangle 
 \end{array}
\right] \, \in \, \R^{2N},
\end{split}
\eeq   
where, for each $i=1,...,N$ and $j=1,2$, the mapping $\G^j_i$ is defined as
\beq 
 \label{5.8}
\left\{\  
\begin{split}
& \Big\langle \G^1_i(\vec u, \vec v, \xi )  ,\, \phi_i \Big\rangle\, :=
\\
& \ \ \ \ \ \int_\Om\Big[ \A_\xi :(\D u_i^\top \D \phi_i) \, +  \big(\K_\xi u_i-S_i\big)\cdot \phi_i \Big] \, \mathrm{d}\mL^n \,+ \int_{\p\Om}\big[ (\ga u_i-s_i)\cdot \phi_i\big] \, \mathrm{d}\mH^{n-1}
\end{split}
\right.
\eeq 
and  
\beq 
 \label{5.9}
\left\{\  
\begin{split}
& \Big\langle \G^2_i(\vec u, \vec v, \xi )  ,\, \psi_i \Big\rangle\, :=
\\
& \ \ \ \ \  \int_\Om\Big[ \A_\xi :(\D v_i^\top \D \psi_i) \, +  \big(\K_\xi v_i- \xi \H u_i\big)\cdot \psi_i \Big] \, \mathrm{d}\mL^n \,+ \int_{\p\Om}\big[ \ga v_i\cdot \psi_i\big] \, \mathrm{d}\mH^{n-1},
\end{split}
\right.
\eeq   
for any test functions 
\[
(\phi_i, \psi_i ) \ \in W^{1,\frac{m}{m-2}} (\Om;\R^{2}) \by W^{1,\frac{p}{p-1}} (\Om;\R^{2}). 
\]
Let us also define for a fixed $M>0$ the next convex weakly closed subset of the Banach space $\mW^p(\Om)$:
\beq
 \label{5.10}
\mW_M^p(\Om) \,:=\, W^{1,\frac{m}{2}}(\Om;\R^{2\by N}) \by W^{1,p}(\Om;\R^{2\by N}) \by W^{2,m}(\Om;[0,M]).
\eeq 
Then, in view of \eqref{5.6}-\eqref{5.10}, we may reformulate the admissible class $\mX^p(\Om)$ of the minimisation problem \eqref{4.15A} as
\beq
\label{5.11}
\mX^p(\Om) \, =\, \Big\{ \big(\vec u, \vec v,\xi\big) \in \mW_M^p(\Om)\, : \ \G \big(\vec u, \vec v,\xi\big)=0\Big\}.
\eeq

With the aim of deriving the variational inequality which is the necessary condition of the minimisation problem \eqref{4.15A}, we compute the Frech\'et derivative of the mapping $\G$ above and prove that it is a submersion.

\begin{lemma} \label{le8}
The mapping $\G$ defined in \eqref{5.6}-\eqref{5.10} is a continuously Frech\'et differentiable submersion and its derivative
\beq \label{5.12}
\begin{split}
\mathrm d \hspace{0.5pt} \G  \ : \ \ \mW^p(\Om) \larrow \mathscr{L} & \left(\mW^p(\Om), \Big[\big( W^{1,\frac{m}{m-2}} (\Om;\R^{2})\big)^*  \by \big(W^{1,\frac{p}{p-1}} (\Om;\R^{2})\big)^*\Big]^N \right)
 \end{split}
\eeq
which maps
\[
(\vec u, \vec v,\xi) \mapsto \big(\mathrm d \hspace{0.5pt} \G\big)_{ (\vec u, \vec v,\xi)}
\]
is given by
\beq  \label{5.13}
\begin{split}
& \left\langle \big(\mathrm d \hspace{0.5pt} \G\big)_{ (\vec u, \vec v,\xi)} (\vec z, \vec w,\eta) ,\, (\vec \phi, \vec \psi\, )\right\rangle\, = \, \left[
\begin{array}{c}
 \Big\langle \big(\mathrm d \hspace{0.5pt} \G^1_1\big)_{ (\vec u, \vec v,\xi)} (\vec z, \vec w,\eta)  ,\,  \phi_1 \Big\rangle \ms
 \\
  \Big\langle \big(\mathrm d \hspace{0.5pt} \G^2_1\big)_{ (\vec u, \vec v,\xi)} (\vec z, \vec w,\eta)  ,\,  \psi_1  \Big\rangle 
 \\
 \vdots
 \\
 \Big\langle \big(\mathrm d \hspace{0.5pt} \G^1_N\big)_{ (\vec u, \vec v,\xi)} (\vec z, \vec w,\eta)   ,\,  \phi_N \Big\rangle \ms
 \\
 \Big\langle \big(\mathrm d \hspace{0.5pt} \G^2_N\big)_{ (\vec u, \vec v,\xi)} (\vec z, \vec w,\eta)  ,\, \psi_N \Big\rangle 
 \end{array}
\right]
\end{split}
\eeq   
where, for each $i \in\{1,...,N\}$ and $j\in\{1,2\}$, we have
\beq 
 \label{5.14}
\left\{\ \
\begin{split}
& \left\langle \big(\mathrm d \hspace{0.5pt} {\G^1_i}\big)_{ (\vec u, \vec v,\xi)} (\vec z, \vec w,\eta) ,\,  \phi_i \right\rangle\, =
\\
\ \ &    \int_\Om\Big[ \A_\xi :(\D z_i^\top \D \phi_i) \, + (\K_\xi z_i)\cdot \phi_i \Big] \, \mathrm{d}\mL^n 
\\
& \ + \int_{\p\Om}(\ga z_i)\cdot \phi_i \, \mathrm{d}\mH^{n-1} \, + \int_\Om \Big[\dot r(\cdot, \xi)\big(  \D u_{i} : \D \phi_{i}\big)  +  u_{i} \cdot \phi_{i}  \Big] \eta \, \mathrm{d}\mL^n
\end{split}
\right.
\eeq 
and  
\beq 
 \label{5.15}
\left\{ \ \ 
\begin{split}
& \left\langle \big(\mathrm d \hspace{0.5pt} {\G^2_i} \big)_{ (\vec u, \vec v,\xi)} (\vec z, \vec w,\eta)  ,\,  \psi_i  \right\rangle\, =
\\
\ \ &   \int_\Om\Big[ \A_\xi :(\D w_i^\top \D \psi_i) \, +  \big(\K_\xi w_i-\H (\eta u_i +\xi z_i )\big)\cdot \psi_i \Big] \, \mathrm{d}\mL^n 
\\
&\ + \int_{\p\Om} (\ga w_i)\cdot \psi_i \, \mathrm{d}\mH^{n-1}\, + \int_\Om \Big[ \dot r(\cdot, \xi)\big( \D v_{i} :\D \psi_{i} \big)  + v_{i} \cdot \psi_{i} \Big] \eta \, \mathrm{d}\mL^n
,
\end{split}
\right.
\eeq   
for any test functions 
\[
(\vec \phi, \vec \psi\, ) \ \in W^{1,\frac{m}{m-2}} (\Om;\R^{2\by N}) \by W^{1,\frac{p}{p-1}} (\Om;\R^{2\by N})
\]
and any $(\vec u, \vec v, \xi ), (\vec z, \vec w, \eta  ) \in \mW^p(\Om)$. 
\end{lemma}

\bp The mapping $\G$ is at most quadratic in all arguments and also continuously Frech\'et differentiable in the space $\mW^p(\Om)$. The form of the derivative can be computed by using directional Gateaux differentiation
\[
 \big(\mathrm d \hspace{0.5pt} \G\big)_{ (\vec u, \vec v,\xi)}(\vec z, \vec w, \eta ) \, = \, \frac{\mathrm d}{\mathrm d \e}\Big|_{\e=0} \G\Big((\vec u, \vec v, \xi )+\e (\vec z, \vec w, \eta )\Big). 
\]
The exact form of the Gateaux derivative of $\G$ is a simple consequence of the definitions of $\A_\xi,\K_\xi$ and the next computations: 
\[
\begin{split}
\frac{\mathrm d}{\mathrm d \e}\Big|_{\e=0}\H (\xi+\e\eta) (u_i+\e z_i) \, &=\, \H (\eta u_i +\xi z_i ),
\\
\frac{\mathrm d}{\mathrm d \e}\Big|_{\e=0}\A_{\xi+\e \eta} : \big((\D u_i+\e\D z_i)^\top \D \phi_i\big) \, &=\, \A : \frac{\mathrm d}{\mathrm d \e}\Big|_{\e=0}(\D u_i+\e\D z_i)^\top \D \phi_i 
\\
&\ \ \ \ + \, \frac{\mathrm d}{\mathrm d \e}\Big|_{\e=0} r\big(\cdot,\xi+\e\eta\big) (\D u_i+\e\D z_i) : \D \phi_i
\\
&=\, \A : (\D z_i^\top \D \phi_i) \, + \, r(\cdot,\xi\big) (\D z_i : \D \phi_i) \phantom{^\big|}
\\
&\ \ \ \, + \, \dot r(\cdot,\xi\big) (\D u_i : \D \phi_i)\, \eta
\\
&=\, \A_\xi : (\D z_i^\top \D \phi_i) \, + \,  \dot r(\cdot,\xi\big) (\D u_i : \D \phi_i)\, \eta \phantom{^\big|},
\end{split}
\]
and
\[
\begin{split}
\frac{\mathrm d}{\mathrm d \e}\Big|_{\e=0}\K_{\xi+\e \eta} (u_i+\e z_i)\cdot \phi_i \, &=\, \K \frac{\mathrm d}{\mathrm d \e}\Big|_{\e=0}( u_i+\e  z_i) \cdot \phi_i 
\\
&\ \ \ \ + \, \frac{\mathrm d}{\mathrm d \e}\Big|_{\e=0} (\xi+\e\eta\big) ( u_i+\e z_i)\cdot \phi_i
\\
&=\, (\K z_i) \cdot \phi_i \, + \, \xi ( z_i \cdot \phi_i)
\, + \, (u_i \cdot \phi_i) \, \eta \phantom{^\big|}
\\
&=\, (\K_\xi z_i) \cdot \phi_i \, + \, (u_i \cdot \phi_i) \, \eta \phantom{^\big|}.
\end{split}
\]
To conclude, we need to show that $\G$ is a submersion, namely that for any $(\vec u, \vec v, \xi ) \in \mW^p(\Om)$, the differential at this point, which is
\[
\big(\mathrm d \hspace{0.5pt} \G\big)_{ (\vec u, \vec v,\xi)} \ : \ \ \mW^p(\Om) \larrow \Big[\big( W^{1,\frac{m}{m-2}} (\Om;\R^{2})\big)^*  \by \big(W^{1,\frac{p}{p-1}} (\Om;\R^{2})\big)^*\Big]^N
\]
is surjective. To this end, for each $i\in\{1,...,N\}$ fix functionals
\[
(\Phi_i, \Psi_i ) \ \in \big( W^{1,\frac{m}{m-2}} (\Om;\R^{2})\big)^*  \by \big(W^{1,\frac{p}{p-1}} (\Om;\R^{2})\big)^*.
\]
This means that there exist
\[
\left\{ \ \
\begin{split}
& (f_i,F_i) \, \in L^{\frac{m}{2}} (\Om;\R^{2})  \by L^{\frac{m}{2}} (\Om;\R^{2\by n}),
\\
&(g_i,G_i) \, \in L^{p} (\Om;\R^{2})  \by L^{p} (\Om;\R^{2\by n}),
\end{split}
\right.
\]
such that the next representation formulas hold true (see e.g.\ \cite{Ad})
\beq
\label{5.16}
\left\{\ \ 
\begin{split}
\langle \Phi_i, \phi_i \rangle \, &=\, \int_\Om \big(f_i\cdot \phi_i \,+\, F_i :\D \phi_i\big)\, \mathrm d\mL^n, 
\\
 \langle \Psi_i, \psi_i \rangle \, &=\, \int_\Om \big(g_i\cdot \psi_i \,+\, G_i :\D \psi_i\big)\, \mathrm d\mL^n,
 \end{split}
 \right.
\eeq
for any $\phi_i \in W^{1,\frac{m}{m-2}} (\Om;\R^{2})$ and $\psi_i \in W^{1,\frac{p}{p-1}} (\Om;\R^{2})$. Then, by \eqref{5.12}-\eqref{5.16}, the surjectivity of the differential $\G'(\vec u, \vec v, \xi )$ follows from the solvability in $(z_i,w_i)$ (for the choice $\eta \equiv 0$) in the weak sense of the PDE systems
\[
\left\{\ \ 
\begin{array}{ll}
  -\div(\A_\xi  \hspace{1pt} \overset{_\bullet}{\phantom{.}}  \D z_i)\,+\, \K_\xi z_i\, =\, f_i-\div F_i, & \ \ \text{ in }\Om, \ms
\\
 (\A  \hspace{1pt} \overset{_\bullet}{\phantom{.}}  \D z_i - F_i)\hspace{1pt} \mathrm n\,+\, \ga z_i\, =\, 0, & \ \ \text{ on }\p\Om,
\end{array}
\right.
\]
and
\[
\left\{\ \ 
\begin{array}{ll}
  -\div(\A_\xi  \hspace{1pt} \overset{_\bullet}{\phantom{.}}  \D w_i)\,+\, \K_\xi w_i\, =\, \big( g_i+\xi\H z_i \big)-\div \,G_i, & \ \ \text{ in }\Om, \ms
\\
 (\A_\xi  \hspace{1pt} \overset{_\bullet}{\phantom{.}}  \D w_i - G_i)\hspace{1pt}\mathrm n\,+\, \ga w_i\, =\, 0, & \ \ \text{ on }\p\Om,
\end{array}
\right.
\]
for all $i\in\{1,...,N\}$ and with $\A,\K,\H,u_i,\xi,\ga,f_i,F_i,g_i,G_i$ being fixed coefficients and parameters. The solvability of the above systems follows from Theorems \ref{th1}-\ref{th3}. The result is therefore complete.
\ep

Now we derive the variational inequality through the generalised Kuhn-Tucker theory of Lagrange multipliers.

\begin{proposition}[The variational inequality] \label{pr9}
For any $p>2n/(n-2)$, there exist Lagrange multipliers
\[
\big(\vec \phi_p, \vec \psi_p \big) \ \in W^{1,\frac{m}{m-2}} (\Om;\R^{2\by N}) \by W^{1,\frac{p}{p-1}} (\Om;\R^{2\by N})
\]
associated with the constrained minimisation problem \eqref{4.15A} for $\E_p$ in the admissible class \eqref{5.11}, such that the constrained minimiser $\big(\vec u_p, \vec v_p, \xi_p\big) \in \mX^p(\Om)$ satisfies the inequality
\beq
\label{5.17}
\begin{split}
  \frac{1}{p} \big(\mathrm d \hspace{0.5pt} \E_p\big)_{ (\vec u_p, \vec v_p,\xi_p)}   \big(\vec z, \vec w, \eta-\xi_p\big) 
 \, & \geq  \, \sum_{i=1}^N \bigg\langle \big(\mathrm d \hspace{0.5pt} \G^1_i\big)_{ (\vec u_p, \vec v_p,\xi_p)}\big(\vec z, \vec w,\eta-\xi_p\big),\, \phi_{pi}   \bigg\rangle
\\
& + \, \sum_{i=1}^N \bigg\langle \big(\mathrm d \hspace{0.5pt} \G^2_i\big)_{ (\vec u_p, \vec v_p,\xi_p)}\big(\vec z, \vec w,\eta-\xi_p\big),\, \psi_{pi} \bigg\rangle,
\end{split}
\eeq  
for any $(\vec z, \vec w,\eta)$ in the convex set $\mW^p_M(\Om)$ (see \eqref{5.10}).
\end{proposition}

\bp In view of Lemmas \ref{pr7}-\ref{le8}, $\E_p$ is Frech\'et differentiable and $\G$ is a continuously Frech\'et differentiable submersion everywhere on $\mW^p(\Om)$, whilst the set $\mW^p_M(\Om)$ is convex and with non-empty interior (with respect to the norm topology). Hence, the hypotheses of the generalised Kuhn-Tucker theory are satisfied (see e.g.\ \cite[p.\ 417-418, Theorem 48B \& Corollary 48.10]{Z}). Therefore, there exists a Lagrange multiplier
\[
\begin{split}
\Lambda_p \ \in \ & \Big(\big(W^{1,\frac{m}{m-2}} (\Om;\R^{2\by N})\big)^* \by \big(W^{1,\frac{p}{p-1}} (\Om;\R^{2\by N})\big)^*\Big)^*
\end{split}
\]
which by standard duality arguments regarding product Banach spaces and their dual spaces that it can be identified with a pair of functions
\[
\big(\vec \phi_p, \vec \psi_p \big) \ \in W^{1,\frac{m}{m-2}} (\Om;\R^{2\by N}) \by W^{1,\frac{p}{p-1}} (\Om;\R^{2\by N})
\]
such that, the constrained minimiser $\big(\vec u_p, \vec v_p, \xi_p\big)$ satisfies
\beq \label{5.19}
\begin{split}
  \frac{1}{p}  \big(\mathrm d \hspace{0.5pt} \E_p   &  \big)_{ (\vec u_p, \vec v_p,\xi_p)}    \Big(\vec z -\vec u_p, \, \vec w -\vec v_p,\, \eta-\xi_p\Big)
 \\
 & - \, \sum_{i=1}^N \bigg\langle \big(\mathrm d \hspace{0.5pt} \G^1_i\big)_{ (\vec u_p, \vec v_p,\xi_p)} \Big(\vec z -\vec u_p, \, \vec w -\vec v_p,\, \eta-\xi_p\Big),\, \phi_{pi}   \bigg\rangle
\\
& - \, \sum_{i=1}^N \bigg\langle \big(\mathrm d \hspace{0.5pt} \G^2_i\big)_{ (\vec u_p, \vec v_p,\xi_p)} \Big(\vec z -\vec u_p, \, \vec w -\vec v_p,\, \eta-\xi_p\Big),\, \psi_{pi} \bigg\rangle\, \geq\, 0,
\end{split}
\eeq
for any $(\vec z, \vec w,\eta)$ in the convex set $\mW^p_M(\Om)$. Recall now that \eqref{5.10} implies that the convex subset $\mW^p_M(\Om)$ of the Banach space $\mW^p(\Om)$ can be written as the Cartesian product of the vector spaces 
\[
W^{1,\frac{m}{2}} (\Om;\R^{2\by N}) \by W^{1,p} (\Om;\R^{2\by N})
\]
with the convex set $W^{2,m}(\Om,[0,M]) $, we may replace $\vec z$ by $\vec z +\vec u_p$ and we may also replace $\vec w$ by $\vec w +\vec v_p$ in \eqref{5.19} to arrive at \eqref{5.17}. The proof of Proposition \ref{pr9} is now complete. 
\ep
We now use Proposition \ref{pr9} to deduce that the variational inequality takes the form \eqref{5.18} below, as a direct consequence of Lemmas \ref{pr7}-\ref{le8}, \eqref{5.1}, \eqref{5.3}, \eqref{5.2}, \eqref{5.12}-\eqref{5.15}. 

\begin{corollary} \label{cor10}
In the setting of Proposition \ref{pr9}, in view of the form of the Frech\'et derivatives of $\E_p$ and $\G$, the variational inequality \eqref{5.17} takes the form
\beq
\label{5.18}
\begin{split}
& \int_{\p\Om} \vec w : \mathrm d [\vec \nu_p(\vec v_p)]\, + \, \frac{\al m}{p} \int_\Om \big(\D^2\eta-\D^2\xi_p \big) : \mu(\D^2\xi_p)\, \mathrm d \mL^n
\\
& \geq \, \sum_{i=1}^N \bigg\{\int_\Om\Big[ \A :(\D z_i^\top \D \phi_{pi}) \, + (\K z_i)\cdot \phi_{pi} \Big] \, \mathrm{d}\mL^n \,+ \int_{\p\Om}(\ga z_i)\cdot \phi_{pi} \, \mathrm{d}\mH^{n-1}\bigg\}
\\
& \ \ \ + \sum_{i=1}^N \bigg\{\int_\Om \Big[ \B :(\D w_i^\top \D \psi_{pi}) \, +  \Big(\L w_i- \H \big((\eta-\xi_p) u_{pi} +\xi_p z_i \big)\Big) \cdot \psi_{pi} \Big] \, \mathrm{d}\mL^n 
\\
&\ \ \ + \int_{\p\Om} (\ga w_i)\cdot \psi_{pi} \, \mathrm{d}\mH^{n-1} \bigg\} +\, \sum_{i=1}^N  \int_\Om (\eta- \xi_p) \bigg(\dot r(\cdot, \xi_p)\Big[  \D u_{pi} : \D \phi_{pi}
\\
&\ \ \ + \, \D v_{pi} :\D \psi_{pi} \Big] +  u_{pi} \cdot \phi_{pi} + v_{pi} \cdot \psi_{pi}    \bigg) \mathrm{d}\mL^n
\end{split}
\eeq
for any $ (\vec z, \vec w, \eta  ) \in \mW^p_M(\Om)$. 
\end{corollary}

We conclude this section by obtaining the further desired information on the variational inequality \eqref{5.18}.

\begin{lemma} \label{le11}
In the setting of Corollary \ref{cor10}, the variational inequality \eqref{5.18} for the constrained minimiser $\big(\vec u_p, \vec v_p ,\xi_p \big)$ is equivalent to the triplet of relations \eqref{5.20}-\eqref{5.22}.
\end{lemma}

\bp The inequality \eqref{5.20} follows by setting $\vec z =\vec w =0$ in \eqref{5.18}. The identity \eqref{5.21} follows by setting $\eta =\xi_p$ and $\vec z =0$ in \eqref{5.18} and by recalling that $W^{1,p} (\Om;\R^{2\by N})$ is a vector space, so the inequality we obtain in fact holds for both $\pm w$. Finally, the identity \eqref{5.22} follows by setting $\eta =\xi_p$ and $\vec w =0$ in \eqref{5.18} and by recalling again that $W^{1,\frac{m}{2}} (\Om;\R^{2\by N})$ is a vector space, so the inequality holds for both $\pm z$.
\ep

\ms

\section{Kuhn-Tucker theory and Lagrange multipliers for the $\infty$-error}  \label{Section6}

In this section we consider the $L^\infty$-minimisation problem \eqref{4.15B} solved in part (B) of Theorem \ref{th4} (Section \ref{Section4}). The goal is to derive the relevant variational inequalities associated with the constrained minimisation of the functional $\E_\infty$ (see \eqref{4.5}) in the admissible class \eqref{4.9}, by analogy to the results in Theorem \ref{th7} of Section \ref{Section5}. To this aim, let us set
\[
C_\infty:=\, \limsup_{p_j \to \infty}C_p,
\]
where
\[
 C_p\,:=\, \|\vec \phi_p\|_{W^{1,\frac{m}{m-2}} (\Om)} +\,  \| \vec \psi_p \|_{W^{1,1}(\Om)}
\]
and $\big(\vec \phi_p, \vec \psi_p \big)$ are the Lagrange multipliers associated with the constrained minimisation problem \eqref{4.15B} (Theorem \ref{th7}). The main result here is therefore the following.

\begin{theorem}[The variational inequalities in $L^\infty$] \label{th14}

In the setting of Section \ref{Section5} and under the same assumptions, suppose additionally that $m>2n$ and also
\beq
\label{6.1a}
\A \in C^1\big(\Om;\R^{n\by n}_+\big), \ \ \ \K,\H \in W^{1,\infty}\big(\Om;\R^{2\by 2}\big), \ \ \ \ka \in C^1(\Om), \ \ \ \vec S \in W^{1,\frac{m}{2}}(\Om;\R^{2\by N}).
\eeq
Then, there exists a subsequence $(p_j)_1^\infty$ and a limiting measure
\[
\vec \nu_\infty \, \in  \mM(\p\Om;\R^{2\by N})
\]
such that
\beq
\label{6.2}
 \vec \nu_p(\vec v_p)  \, \weakstar \, \vec \nu_\infty \ \text{ in }\mM(\p\Om;\R^{2\by N}),
\eeq
as $p_j \to \infty$, where $\vec \nu_p(\vec v_p)$ is given by \eqref{5.2}, and:

\ms

\noi {\rm (I)} If $C_\infty=0$, then $\vec \nu_\infty = \vec 0$.

\ms

\noi {\rm (II)} If $C_\infty>0$, then there exist (rescaled) limiting Lagrange multipliers
\[
\big(\vec \phi_\infty, \vec \psi_\infty \big) \ \in W^{1,\frac{m}{m-2}} (\Om;\R^{2\by N}) \by BV (\Om;\R^{2\by N}),
\]
such that
\beq
\label{6.1}
 \bigg( \frac{\vec \phi_p}{C_p} \, , \frac{\vec \psi_p}{C_p} \bigg) \ \weakstar \ \big(\vec \phi_\infty, \vec \psi_\infty \big) \ \ \text{ in }\ W^{1,\frac{m}{m-2}} (\Om;\R^{2\by N}) \by BV (\Om;\R^{2\by N})
\eeq
as $p_j \to \infty$. Then, for the above Lagrange multipliers, the constrained minimiser $\big(\vec u_\infty, \vec v_\infty, \xi_\infty\big) \in \mX^\infty(\Om)$ satisfies the next three relations:
\beq
\label{6.3}
\begin{split}
 \sum_{i=1}^N    \int_\Om  (\eta- \xi_\infty) & \bigg[ \big(\H  u_{\infty i}  \big)  \cdot \psi_{\infty i} -\Big( u_{\infty i} \cdot \phi_{\infty i} + v_{\infty i} \cdot \psi_{\infty i} \Big)  \bigg)\,  \mathrm{d}\mL^n 
\\
& \geq \,   \sum_{i=1}^N \int_\Om  (\eta- \xi_\infty)  \dot r(\cdot, \xi_\infty) 
\Big( \D u_{\infty i}: \D \phi_{\infty i}\Big) \mathrm{d}\mL^n  
\\
& \ \ \ + \,  \sum_{i=1}^N  \int_\Om  (\eta- \xi_\infty)  \dot r(\cdot, \xi_\infty)\, \D v_{\infty i}: \, \mathrm{d}[\D \psi_{\infty i}] 
\end{split}
\eeq
for any $\eta \in C_{\xi_\infty}(\Om;[0,M])$ (namely $\eta \in C(\Om;[0,M])$ with $\eta =\xi_\infty$ on $\p\Om$); further,
\beq
\label{6.4}
\begin{split}
 \frac{1}{C_\infty}\int_{\p\Om} \vec w : \mathrm d \vec \nu_\infty \, =\, & \sum_{i=1}^N  \bigg\{  \int_\Om \A_{\xi_\infty} : (\D w_i)^\top  \mathrm{d}[\D \psi_{\infty i}] \
\\
& +  \int_{\p\Om}  \big(\K_{\xi_\infty} w_i\big) \cdot \psi_{\infty i}  \, \mathrm{d}\mL^n \, + \int_{\p\Om} (\ga w_i)\cdot \psi_{\infty i} \, \mathrm{d}\mH^{n-1} \bigg\},
\end{split}
\eeq
for any $\vec w  \in C_0^1(\overline{\Om};\R^{2\by N})$, and finally
\beq
\label{6.5}
\begin{split}
  \sum_{i=1}^N \int_\Om\Big[ \A_{\xi_\infty} :(\D z_i^\top \D \phi_{\infty i}) \,  +  (\K_{\xi_\infty}  &  z_i)  \cdot \phi_{\infty i} - \xi_\infty \big(\H z_i \big) \cdot \psi_{\infty i} \Big] \, \mathrm{d}\mL^n 
\\
& =\, -\, \sum_{i=1}^N \int_{\p\Om}(\ga z_i)\cdot \phi_{\infty i} \, \mathrm{d}\mH^{n-1} 
\end{split}
\eeq
for any $ \vec z  \in C^1(\overline{\Om};\R^{2\by N})$. 
\end{theorem}

Note the interesting fact that the limiting variational inequality \eqref{6.3} has no dependence on the regularisation parameter $\al$, as the corresponding term in \eqref{5.20} is annihilated.

\bp We begin by showing that for any $p>n$ and any $\vec v  \in W^{1,p}(\Om;\R^{2\by N})$, we have the next total variation bound for the measure \eqref{5.2}:
\begin{align}
\big\| \vec \nu_p(\vec v)\big\|(\p\Om)\, &\leq\, N. \label{5.4}
\end{align}
Indeed, by H\"older inequality we have
\[
\begin{split}
\big\| \nu_{p_i}(\vec v) \big\|(\p\Om) \, & \leq\, \frac{ \av_{\Ga_i}\big(|v_i-v^\de_i|_{(p)}\big)^{p-2} |v_i-v^\de_i| \, \mathrm d  \mH^{n-1}  }{\left(\, \av_{\Ga_i}\big(\big| v_i -v^\de_i \big|_{(p)} \big)^p\, \mathrm d \mH^{n-1}\right)^{\!\!\frac{p-1}{p}}}
\\
&  \leq\, \frac{ \av_{\Ga_i}\big(|v_i-v^\de_i|_{(p)}\big)^{p-1} \, \mathrm d  \mH^{n-1}  }{\left(\, \av_{\Ga_i}\big(\big| v_i -v^\de_i \big|_{(p)} \big)^p\, \mathrm d \mH^{n-1}\right)^{\!\!\frac{p-1}{p}}}
\\
&\leq \, 1,
\end{split}
\]
for any $i\in\{1,...,N\}$. By the sequential weak* compactness of the spaces of Radon measures, the estimate \eqref{5.4} implies the existence of a subsequence $(p_j)_1^\infty$ and of the claimed limit measure $\vec \nu_\infty$ in \eqref{6.2}.

Now we proceed with establishing (I) and (II) of the theorem.

\smallskip

\noi (I) Suppose that $C_\infty=0$. Then, it follows that 
\beq
\label{6.9}
\ \ \  \big( \vec \phi_p, \vec \psi_p\big) \larrow \big(\vec 0, \vec 0 \big)\ \ \text{ in }  W^{1,\frac{m}{m-2}} (\Om;\R^{2\by N})\! \by BV (\Om;\R^{2\by N})
\eeq
as $p_j \to \infty$, where $\big(\vec \phi_p, \vec \psi_p \big)$ are the Lagrange multipliers associated with the constrained minimisation problem \eqref{4.15B}. By \eqref{6.2} and \eqref{6.9}, by passing to the limit as $p_j \to \infty$ in \eqref{5.21}, we obtain
\[
\begin{split}
 \int_{\p\Om} \vec w : \mathrm d \vec \nu_\infty \, =\,0 \, =\, \lim_{p_j\to \infty} \sum_{i=1}^N & \bigg\{  \int_\Om \Big[ \A_{\xi_p} :(\D w_i^\top \D \psi_{pi}) \, +  \big(\K_{\xi_p} w_i\big) \cdot \psi_{pi} \Big] \, \mathrm{d}\mL^n 
\\
& \ + \int_{\p\Om} (\ga w_i)\cdot \psi_{pi} \, \mathrm{d}\mH^{n-1} \bigg\}, 
\end{split}
\]
for any $\vec w  \in C^1 (\overline{\Om};\R^{2\by N})$. Therefore, $\vec \nu_\infty =\vec 0$, as claimed.

\ms

\noi (II) Suppose now $C_\infty>0$. Then, the desired relations \eqref{6.3}-\eqref{6.5} follow directly from \eqref{5.21}-\eqref{5.22} by rescaling the Lagrange multipliers $\big(\vec \phi_p ,\vec \psi_p \big)$ and passing to the limit as $p_j \to \infty$, since the rescaled Lagrange multipliers $\big(\vec \phi_p/C_p ,\vec \psi_p/C_p \big)$ are bounded in the product space
\[
W^{1,\frac{m}{m-2}} (\Om;\R^{2\by N}) \by BV (\Om;\R^{2\by N})
\]
and therefore the sequence is weakly* precompact. (Recall also that on a reflexive space the weak and the weak* topology coincide.) Note first that we have
\[
\begin{split}
\bigg|\! \int_\Om \D^2(\eta - \xi_p): \mu\big(\D^2\xi_p\big) \, \mathrm d \mL^n \bigg|\, & =\, \Bigg|\, {\av_\Om} \D^2(\eta - \xi_p\big): \frac{ (|\D^2\xi_p|_{(m)})^{m-2} \, \D^2\xi_p}{\big(\|\D^2\xi_p\|_{\dot L^m(\Om)}\big)^{m-1}} \, \mathrm d \mL^n \Bigg|
\\
& \leq \, {\av_\Om} \big|\D^2(\eta - \xi_p)\big|  \frac{ (|\D^2\xi_p|_{(m)})^{m-1}}{\big(\|\D^2\xi_p\|_{\dot L^m(\Om)}\big)^{m-1}} \, \mathrm d \mL^n 
\\
& \leq \, C \big\|\D^2(\eta - \xi_p)\big\|_{L^m(\Om)} ,
\end{split}
\]
by the definition of $\mu$ and by H\"older inequality. In order to conclude, we need to justify the weak* convergence as $p_j \to \infty$ of the quadratic terms 
\[
\D u_{pi} :\frac{\D\phi_{pi}}{C_p}\ , \ \ \ \D v_{pi}:\frac{\D\psi_{pi}}{C_p}.
\]
To this end, we will show that under the higher regularity assumptions on the coefficients, we in fact have the next strong modes of convergence for the $p$-minimisers:
\begin{align}
\label{6.8} \D u_{pi} & \larrow \D u_{\infty i} \ \  \text{ in }L_{\text{loc}}^{\frac{m}{2}}(\Om;\R^2),
\\
\label{6.9} \D v_{pi} & \larrow \D v_{\infty i} \ \ \, \text{ in }C(\Om;\R^2),
\end{align}
as $p_j \to \infty$, for all $i\in\{1,...,N\}$. Before proving \eqref{6.8}-\eqref{6.9}, we demonstrate how to conclude by assuming them. Since we have
\begin{align}
\label{6.10} \frac{ \D\phi_{p i}}{C_p}  & \weak \,  \D\phi_{\infty i} \ \ \text{ in }L^{\frac{m}{m-2}}(\Om;\R^2),
\\
\label{6.11} \frac{ \D\psi_{p i}}{C_p}\mL^n  & \,  \weakstar \,  \D\psi_{\infty i} \ \ \ \text{ in }\mM(\Om;\R^2)
\end{align}
and also $\xi_p \larrow \xi_\infty$ in $C^1(\overline{\Om})$ as $p_j \to \infty$, by choosing any $\mO\Subset \Om$ with Lipschitz boundary (for instance the union of finitely many balls), $(\eta_{p_j})_1^\infty \sub W^{2,m}(\Om;[0,M])$ and $\eta \in W^{2,m}(\Om;[0,M])$ with
\[
\eta_p \equiv \xi_p \ \text{ on }\Om\set\mO, \ \ \ \eta_p \larrow \eta \, \text{ in }\, W^{2,m}(\mO) \sub C^1(\overline{\mO}),
\]
we have $\eta -\xi_p \in W^{2,m}_0(\mO)$ and 
\[
\eta_p-\xi_p \larrow \eta-\xi_\infty \, \text{ in }\, W^{2,m}_0(\mO)
\]
as $p_j \to \infty$. Hence, \eqref{6.3}-\eqref{6.5} follow by \eqref{6.8}-\eqref{6.11}, together with the weak-strong continuity of the duality pairing between $L^{\frac{m}{2}}(\mO;\R^2)$ and $L^{\frac{m}{m-2}}(\mO;\R^2)$ and the weak*-strong continuity of the duality pairing between $C_0(\mO;\R^2)$ and $\mM(\mO;\R^2)$, at least for test functions $\eta \in W^{2,m}_{\xi_\infty}(\mO;[0,M])$. The general case for test functions $\eta \in C_{\xi_\infty}(\Om;[0,M])$ follows by a standard approximation argument.

Now we establish \eqref{6.8}-\eqref{6.9}. Fix $i\in\{1,...,N\}$, $e\in \R^n$ with $|e|=1$,  $h\neq0$ and $\big( \vec u, \vec v, \xi )\in \mX^p(\Om)$ for some $p$ large. Fix also $\zeta \in C^1_c(\Om)$ and let $\D^{1,h}_e$ symbolise the difference quotient with step size $h$ along the direction of $e$. By testing in the weak form of the equations $(a)_i-(b)_i$ appearing in the constrained class \eqref{4.8} against test functions of the form
\[
\phi_i \,=\, \psi_i \,:=\, \D^{1,-h}_e \zeta ,
\]
standard regularity arguments imply that the directional derivatives $\D_e u_i$ and $\D_e v_i$ solve weakly the divergence systems
\begin{align}
\nonumber \ \  -\div(\A_\xi  \hspace{1pt} \overset{_\bullet}{\phantom{.}}  \D (\D_eu_i))\,=& \,   \Big\{\! - \K_\xi (\D_eu_i)+ \D_e S_i  - (\D_e\K_\xi) u_i\Big\} +\div\big( (\D_e \A_\xi)  \hspace{1pt} \overset{_\bullet}{\phantom{.}}  \D u_i\big) ,
\\
&  \label{6.13} \hspace{220pt} \text{ in }\Om,  
\end{align}
and
\begin{align}
\ \     -\div(\A_\xi \hspace{1pt} \overset{_\bullet}{\phantom{.}}  \D (\D_e v_i))\,= &\,  \Big\{ \! - \K_\xi (\D_e v_i)+ \xi \H (\D_eu_i) + \D_e(\xi \H) u_i -  (\D_e\K_\xi) v_i \Big\} \nonumber
\\
& \label{6.14} + \div\big( (\D_e \A_\xi)  \hspace{1pt} \overset{_\bullet}{\phantom{.}}  \D v_i\big) ,  \hspace{130pt} \text{ in }\Om.
\end{align}
In view of \eqref{1.4}-\eqref{1.5} we have
\[
\D_e \A_\xi \, =\, \D_e \A \,+\, \D_e(r(\cdot,\xi))\mathrm I_n, \ \ \ \ \D_e \K_\xi \, =\, \D_e \K \,+\, (\D_e\xi) \mathrm I_2,
\]
and
\[
\D_e(r(\cdot,\xi))  \, =\, -\frac{\la \big(\D_e\ka  +\D_e\xi \big)}{(\ka +\xi )^2}.
\]
Due to our assumption \eqref{6.1a} and the embedding $W^{2,m}(\Om) \sub C^1(\overline{\Om})$, we have that 
\[
\D_e \A_\xi \, \in C\big(\Om;\R^{n\by n}_+ \big), \ \ \ \D_e \K_\xi \in L^\infty(\Om;\R^{2\by 2} \big)
\]
as
\[
\text{$\D_e(r(\cdot,\xi)) \in C(\Om)$, \ $\K \in W^{1,\infty}(\Om;\R^{2\by 2} \big)$, \  $\xi,\ka \in C^1(\Om)$.}
\]
Further, in view of our assumptions and H\"older inequality we have that 
\[
\left.
\begin{array}{r}
 - \K_\xi (\D_eu_i)+ \D_e S_i  - (\D_e\K_\xi) u_i
 \ms
 \\
 (\D_e \A_\xi)  \hspace{1pt} \overset{_\bullet}{\phantom{.}}  \D u_i
 \\
 \ms
 - \K_\xi (\D_e v_i)+ \xi \H (\D_eu_i) + \D_e(\xi \H) u_i -  (\D_e\K_\xi) v_i
 \\
  (\D_e \A_\xi)  \hspace{1pt} \overset{_\bullet}{\phantom{.}}  \D v_i
 \end{array}
 \right\} 
  \, \in \, L^{\frac{m}{2}}(\Om) 
\]
because
\[
\begin{split}
\xi,\, |\H|,\, |\D_e \H|,\, |\D_e(\xi\H)|, \,|\D_e\K_\xi|,\,|\K_\xi|, \, |\D_e \A_\xi| \, & \in L^\infty(\Om),
\phantom{\big|} \\
 |u_i|,\,|\D_e u_i|,|\,v_i|,\,|\D_e v_i|,\, |\D_e S_i|,\, |\D u_i|,\, |\D v_i| \, & \in  L^{\frac{m}{2}}(\Om),
 \end{split}
\]
for $p>2m$. By the interior $L^2$ and $L^{{\frac{m}{2}}}$ regularity estimates for the systems \eqref{6.13}-\eqref{6.14} (see e.g.\ \cite[Sections 4.3.1 \& 7.1.2]{GM}), for any $\mO \Subset \Om$ there exists $C>0$ independent of $p$ such that
\[
\big\| \D^2 \vec u_{p} \big\|_{L^\frac{m}{2}(\mO)} + \, \big\| \D^2 \vec v_{p} \big\|_{L^\frac{m}{2}(\mO)} \,\leq\, C.
\]
Since by assumption $m>2n$, by the Morrey estimate we have
\[
\big\| \D \vec u_{p} \big\|_{C^{0,1-\frac{2n}{m}}(\mO)} + \, \big\| \D \vec v_{p} \big\|_{C^{0,1-\frac{2n}{m}}(\mO)} \,\leq\, C.
\]
By standard compact embedding arguments in H\"older spaces, \eqref{6.8}-\eqref{6.9} ensue as a consequence of the above estimates. The proof of the theorem is complete.
\ep

\subsection*{Acknowledgments} The author would like to thank the anonymous referee for their constructive critique on an earlier version of this paper, which improved both the content and the presentation.

%\ms

\bibliographystyle{amsplain}

\end{document}